\documentclass[]{amsart} \usepackage{amsfonts}

\usepackage{bbding}
\usepackage[dvips]{graphicx}
\usepackage{amsfonts,amssymb,amsmath,amsthm,amscd}

\usepackage{mathrsfs}
\usepackage{xcolor}
\usepackage{empheq}
\usepackage{adforn}
\usepackage{cancel}
\usepackage{xr}
\usepackage{mdframed}
\usepackage{tikz-cd}
\usepackage{tikz}
\usepackage{amsfonts}
\usepackage{mathdots}
\usepackage{pinlabel}

\usepackage{enumerate}
\usepackage{lmodern}
\usepackage{mdframed}
\usepackage{stmaryrd}
\usepackage{blindtext}
\usepackage{leftidx}
\usepackage{easytable}

\usepackage{fouriernc}
\usepackage[scaled=0.875]{helvet}
\usepackage{courier}
\usepackage[marginparwidth=2.5cm]{geometry}
\definecolor{navie}{RGB}{0, 70, 140}
\usepackage[bookmarks=true,colorlinks=true,urlcolor=red!60!black,linkcolor=red!60!black,citecolor=navie,pdfstartview=FitV]{hyperref}

\usepackage{enumitem}

\makeatletter
\newcommand*{\enuma}[1]{%
	\expandafter\@enuma\csname c@#1\endcsname%
}
\newcommand*{\@enuma}[1]{%
	$\ifcase#1\or(a)\or(b_-)\or(b_+)\or(c)%
	\else\@ctrerr\fi$%
}
\AddEnumerateCounter{\enuma}{\@enuma}{53.13}
\makeatother

\makeatletter
\newcommand*{\enumb}[1]{%
	\expandafter\@enumb\csname c@#1\endcsname%
}
\newcommand*{\@enumb}[1]{%
	$\ifcase#1\or(b_-)\or(b_+)%
	\else\@ctrerr\fi$%
}
\AddEnumerateCounter{\enumb}{\@enumb}{53.13}
\makeatother

\usepackage{enumerate}
\newlist{steps}{enumerate}{1}
\setlist[steps, 1]{itemsep=8pt,leftmargin=0cm,itemindent=.5cm,labelwidth=\itemindent,labelsep=0cm,align=left,label = \textbf{\emph{Step \arabic*}:\,}}

\makeatletter
\newtheorem*{rep@theorem}{\rep@title}
\newcommand{\newreptheorem}[2]{%
\newenvironment{rep#1}[1]{%
 \def\rep@title{#2 \ref{##1}}%
 \begin{rep@theorem}}%
 {\end{rep@theorem}}}
\makeatother

\makeatletter
\newtheorem*{rep@cor}{\rep@title}
\newcommand{\newrepcor}[2]{%
\newenvironment{rep#1}[1]{%
 \def\rep@title{#2 \ref{##1}}%
 \begin{rep@cor}}%
 {\end{rep@cor}}}
\makeatother

\makeatletter
\newtheorem*{rep@prop}{\rep@title}
\newcommand{\newrepprop}[2]{%
\newenvironment{rep#1}[1]{%
 \def\rep@title{#2 \ref{##1}}%
 \begin{rep@prop}}%
 {\end{rep@prop}}}
\makeatother

\newtheorem{corollary}{Corollary}[section]
\newtheorem{corx}{Corollary}

\newtheorem{theorem}[corollary]{Theorem}
\newtheorem{thmx}[corx]{Theorem}

\newtheorem{proposition}[corollary]{Proposition}

\newrepcor{cor}{Corollary}
\newreptheorem{theorem}{Theorem}
\newrepprop{prop}{Proposition}

\newtheorem*{theorem*}{Theorem}
\newtheorem{lemma}[corollary]{Lemma}

\newtheorem*{lemma*}{Lemma}

\newtheorem{manualtheoreminner}{Theorem}
\newenvironment{manualtheorem}[1]{%
  \renewcommand\themanualtheoreminner{#1}%
  \manualtheoreminner
}{\endmanualtheoreminner}

\theoremstyle{definition} 
\newtheorem{definition}[corollary]{Definition}
\newtheorem{propdef}[corollary]{Proposition/Definition}

\theoremstyle{remark} \newtheorem{remark}[corollary]{Remark} \numberwithin{equation}{section}
\newtheorem*{remark*}{Remark}
\numberwithin{figure}{section}

\renewcommand{\phi}{\varphi}

\newcommand*\circled[1]{\tikz[baseline=(char.base)]{
		\node[shape=circle,draw,inner sep=0.6pt] (char) {#1};}}
\newcommand{\R}{\mathbb{R}}

\newcommand{\Rp}{\mathbb{R}\cup\{+\infty\}}

\newcommand{\gra}[1]{\mathrm{gr}({#1})}
\newcommand{\sepi}[1]{\mathrm{ep\!}^\circ({#1})}
\newcommand{\dom}[1]{\mathsf{dom}({#1})}
\newcommand{\env}[1]{\overline{#1}}
\newcommand{\interior}{\mathsf{int}\,}

\newcommand{\eps}{\varepsilon}
\newcommand{\C}{\mathsf{C}}

\newcommand{\MA}[1]{\mathcal{M}\left[{#1}\right]}

\newcommand{\ccdot}{\,\cdot\,}

\newcommand{\SC}{\mathsf{S}}

\newcommand{\dif}{\mathsf{d}}

\newcommand{\transp}[1]{\leftidx{^\mathsf{t}}{\!#1}}

\newcommand{\pa}{\partial}

\newcommand{\dx}{{\dif x}}

\newcommand{\eg}{\textit{e.g.\@ }}
\newcommand{\cf}{\textit{c.f.\@ }}
\newcommand{\ie}{\textit{i.e.\@ }}

\newcommand{\D}{\mathsf{D}}

\newcommand{\LC}{\mathsf{LC}}

\newcommand{\A}{\mathbb{A}}

\begin{document}

\title[Hypersurfaces of constant Gauss-Kronecker curvature with Li-normalization]{Hypersurfaces of constant Gauss-Kronecker curvature with Li-normalization in affine space} 
\author{Xin Nie}
\address{Xin Nie: Shing-Tung Yau Center of Southeast University, Nanjing 210018, China.} \email{nie.hsin@gmail.com}
\author{Andrea Seppi}
\address{Andrea Seppi: Institut Fourier, UMR 5582, Laboratoire de Math\'ematiques,
Universit\'e Grenoble Alpes, CS 40700, 38058 Grenoble cedex 9, France.} \email{andrea.seppi@univ-grenoble-alpes.fr}

\thanks{The second author is member of the national research group GNSAGA}

\maketitle



\begin{abstract}
For convex hypersurfaces in the affine space $\mathbb{A}^{n+1}$ ($n\geq2$), A.-M.\ Li introduced the notion of $\alpha$-normal field as a generalization of the affine normal field. By studying a Monge-Amp\`ere equation with gradient blowup boundary condition, we show that regular domains in $\mathbb{A}^{n+1}$, defined with respect to a proper convex cone and satisfying some regularity assumption if $n\geq3$, are 
foliated by complete convex hypersurfaces with constant Gauss-Kronecker curvature relative to the Li-normalization.  When $n=2$, a key feature is that no regularity assumption is required, and the result extends our recent work about the $\alpha=1$ case.
\end{abstract}

\section{Introduction}
As shown by B.\ Chen, A.-M.\ Li, U.\ Simon and G.\ Zhao \cite{lisimoncrelle,lisimongeomded}, Euclidean-complete convex hypersurfaces of constant affine Gauss-Kronecker curvature in the $(n+1)$-dimensional real affine space $\mathbb{A}^{n+1}$ ($n\geq2$) are governed by the Monge-Amp\`ere problem
\begin{equation}\label{eqn_main}
\begin{cases}
\det\D^2 u=(-w_\Omega)^{-(n+2)}\ \text{ in }\Omega,\\
u|_{\pa\Omega}=\phi,\\
\|\D u(x)\|\to+\infty \text{ as $x\in\Omega$ tends to $\pa\Omega$},
\end{cases}
\end{equation}
where $\Omega\subset\R^n$ is a bounded convex domain, $\phi$ is a function on $\pa\Omega$, and $w_\Omega\in\C^0(\overline{\Omega})\cap\C^\infty(\Omega)$ is the unique convex solution, established by Cheng-Yau \cite{chengyau1}, of another Monge-Amp\`ere equation
\begin{equation}\label{eqn_chengyau}
\begin{cases}
\det\D^2 w=(-w)^{-(n+2)}\ \text{ in }\Omega,\\
w|_{\pa\Omega}=0.
\end{cases}
\end{equation}
More specifically, the domain $\Omega$ and the boundary function $\phi$ determine a convex cone $C\subset\R^{n+1}$ and a \emph{$C$-regular domain} $D\subset\A^{n+1}$, respectively, whereas a solution $u$ of Eq.\eqref{eqn_main} corresponds, via Legendre transformation, to a convex hypersurface $\Sigma\subset D$ asymptotic to $\pa D$ with constant affine Gauss-Kronecker curvature, and the last gradient blowup condition in \eqref{eqn_main} is equivalent to the completeness of $\Sigma$. As a well studied particular case, if $C$ is the light cone in the Minkowski space $\R^{n,1}$, so that $\Omega$ is the unit ball and $w_\Omega(x)=-\sqrt{1-|x|^2}$, then $\Sigma$ is a spacelike hypersurface in $\R^{n,1}$ whose Gauss-Kronecker curvature in the classical sense is constant (see \cite{bs, bon_smillie_seppi,choi-treibergs,treibergs}). In \cite{lisimoncrelle,lisimongeomded} and most of the follow-up works, in order to guarantee the solvability of \eqref{eqn_main}, $\pa\Omega$ and $\phi$ are assumed to be $\C^2$ and $\Omega$ to be strictly convex. 

Recently, we investigated the subject again in \cite{nie-seppi2,nie-seppi}, with emphasis on the geometry of regular domains and group actions, through which the link with \emph{higher Teichm\"uller theory} and the theory of \emph{globally hyperbolic spacetimes} is made. Because of this perspective, a key feature of our works is that the regularity assumptions on $\pa\Omega$ and $\phi$ are relaxed, which is necessary for the applications. But this also forces us to restrict to the $n=2$ case, as it is well-known that for $n\geq3$, non-smooth boundary data can result in undesirable Pogorelov-type singular solutions. Bonsante-Fillastre \cite{bonsante-fillastre} studied the geometry of such solutions in globally hyperbolic spacetimes.


\subsection*{Three-dimensional affine space} The purpose of this paper is twofold. The first is to extend the results in \cite{nie-seppi} to convex surfaces in $\A^3$ relative to the \emph{Li-normalization}. Namely, while the above notion of affine  Gauss-Kronecker curvature for a hypersurface $\Sigma\subset\A^{n+1}$ is defined using the usual \emph{affine normal field} $N_1:\Sigma\to\R^{n+1}$, here we replace $N_1$ by a more general transversal vector field $N_\alpha:\Sigma\to\R^{n+1}$ defined by A.-M.\ Li (see \cite{wu-zhao,xiong-yang,xu,xu-li-li}), which depends on a parameter $\alpha\in\R$ and coincides with $N_1$ when $\alpha=1$. For $\alpha\neq0$, the condition that the Gauss-Kronecker curvature of $\Sigma$ relative to $N_\alpha$ is constant turns out to be similar to equations \eqref{eqn_main} and \eqref{eqn_chengyau}, only with the exponent $-(n+2)$ in both equations replaced by $-\frac{n+2}{\alpha}$. The new equations and the problem of prescribing the curvature relative to $N_\alpha$ have been studied by Wu-Zhao \cite{wu-zhao} under the aforementioned $\C^2$ and strict convexity assumptions. We show the following extension in the $n=2$ case, with the weakest possible regularity assumption on $\phi$. Here, for any extended-real-valued function $f$, we let $\dom{f}$ denote the subset of the domain where $f$ is real-valued. 
\begin{thmx}[Simplified version of Theorem \ref{thm_ck2analytic}]\label{thm_1}
	Let $\Omega\subset\R^2$ be a bounded convex domain satisfying the exterior circle condition, and $\phi:\pa\Omega\rightarrow \mathbb{R}\cup\{+\infty\}$ be a lower semicontinuous function such that $\dom{\phi}$ has at least three points. Then for any $\alpha\in(0,1]$, there exists a unique lower semicontinuous convex function $u:\overline{\Omega}\to\Rp$  which is smooth in the interior $U$ of $\dom{u}$ and satisfies
	\begin{equation}\label{eqn_li}
	\begin{cases}
	\det\D^2 u=\left(-w\right)^{-\frac{4}{\alpha}} \text{ in }U:=\interior \dom{u},\\
	u|_{\pa\Omega}=\phi,\\
	\|\D u(x)\|\to+\infty \text{ as $x\in U$ tends to $\pa U$},
	\end{cases}
	\end{equation}
where $w\in\C^0(\overline{\Omega})\cap\C^\infty(\Omega)$ is the unique convex solution to
$$
\begin{cases}
\det\D^2w=(-w)^{-\frac{4}{\alpha}}\ \text{ in }\Omega,\\
w|_{\pa\Omega}=0.
\end{cases}
$$
Moreover, this $u$ has the property that $\dom{u}$ coincides with the convex hull of $\dom{\phi}$ in $\R^2$, and $u$ coincides with the convex envelope function $\env{\phi}$ of $\phi$ on the boundary of $\dom{u}$.
\end{thmx}
\begin{remark*}[about the assumptions] 
	It has been observed in \cite[Prop.\@ E]{nie-seppi} that if $\Omega$ is not strictly convex, then the gradient blowup condition in \eqref{eqn_li} might not be fulfillable for certain $\phi$ that takes finite values exactly at three points. The \emph{exterior circle condition} on $\Omega$ (\ie the condition that for every $p\in \pa\Omega$, there is a circle passing through $p$ surrounding $\Omega$) is intended as a sufficient condition to guarantee the solvability of \eqref{eqn_li}. Meanwhile, we also show in Proposition \ref{prop_circle} that for any $\alpha>1$, \eqref{eqn_li} is not solvable when $\Omega$ is the unit disk and $\phi$ takes finite values exactly at three points, so the assumption $\alpha\leq1$ is sharp.
\end{remark*}

In order to give an affine-differential-geometric interpretation of the theorem, we first pointed out that although the affine normal field $N_1$ of $\Sigma$ is uniquely determined once a translation-invariant volume form on $\A^{n+1}$ is chosen, the definition of Li's $\alpha$-normal field $N_\alpha$ with $\alpha\neq1$ depends furthermore on the choice of a \emph{vertical unit vector} $v$ in the underlying vector space $\R^{n+1}$ (\cf Remark \ref{remark_covariance} below). In fact, for $\alpha=0$, $N_0$ is exactly the constant vector field given by $v$. With this in mind, we can state the geometric counterpart of Theorem \ref{thm_1} as follows:
\begin{manualtheorem}{\ref*{thm_1}'}[Simplified version of Theorem \ref{thm_ck2}]\label{thm_1'}
	Let $\alpha\in(0,1]$ be a constant, $v\in\R^3$ be a non-zero vector and $C\subset\mathbb{R}^3$ be a proper convex cone containing $v$ such that a planar section $\Omega$ of the dual cone  $C^*\subset\R^{3*}$ satisfies the exterior circle condition. Then every proper $C$-regular domain $D\subset\A^3$ is foliated by smooth, complete, locally strongly convex surfaces asymptotic to $\pa D$, whose Gaussian curvatures with respect to Li's $\alpha$-normal fields (defined using $v$ as the vertical unit vector) are constants ranging from $0$ to $+\infty$. Moreover, the function on $D$ which assigns to each leaf its curvature is $\log$-convex.
\end{manualtheorem}

Theorems \ref{thm_1} and \ref{thm_2} generalize the results in \cite{nie-seppi} about the $\alpha=1$ case.

\subsection*{Higher dimensions} The other purpose of this paper is to give a careful exposition of affine hyperspheres and constant Gauss-Kronecker curvature hypersurfaces relative to Li-normalization, including their link with Monge-Amp\`ere equations. In the process, we improve and simplify the past works \cite{lisimoncrelle, lisimongeomded,wu-zhao,xiong-yang} on the PDE side and also provide a clearer geometric picture. In particular, we prove statements similar to Theorems \ref{thm_1} and \ref{thm_1'} which hold for all dimensions $n\geq2$ (with stronger assumptions required). 
We first state the analytic version:

\begin{thmx}[Simplified version of Theorem \ref{thm_ck1analytic}]\label{thm_2}
For any bounded convex domain $\Omega\subset\R^n$ ($n\geq2$), suppose $\phi\in\C^0(\pa\Omega)$ satisfies the following condition: for every $p\in\pa\Omega$, there exists an affine function $a:\R^n\to\R$ such that $a(p)=\phi(p)$ and $a\geq\phi$ on $\pa\Omega$. Then for any $\gamma>n$ and $\lambda>0$,  there exists a unique convex generalized solution $u$ to the  problem
\begin{equation}
\begin{cases}
\det\D^2u=\lambda(-w)^{-\gamma}\ \text{ in }\Omega\\
u|_{\pa\Omega}=\phi
\end{cases}
\end{equation}
which is in $\C^0(\overline{\Omega})\cap\C^\infty(\Omega)$ and has the gradient blowup property on $\pa\Omega$, where $w$ is the unique convex solution to \eqref{eqn_introgamma} (see Theorem \ref{thm_3} below).
\end{thmx}

The affine-differential-geometric counterpart of Theorem \ref{thm_2} is contained the following statement (the parameters in the two statements are related by $\gamma=\frac{n+2}{\alpha}$):

\begin{manualtheorem}{\ref*{thm_2}'}[Simplified version of Theorem \ref{thm_ck1}]\label{thm_2'}
Let $\alpha\in(0,1+\frac{2}{n})$, $v\in\R^{n+1}$ be a non-zero vector, $C\subset\R^{n+1}$ be a proper convex cone  containing $v$, and $D\subset\A^{n+1}$ be a $C$-regular domain satisfying the following condition: for every $C$-null subspace $L\subset\R^{n+1}$, there exists a translation $C'\subset\A^{n+1}$ of the cone $C$ such that $C'$ contains $D$ and the supporting hyperplane $L'$ of $C'$ parallel to $L$ is an asymptotic hyperplane
 of $D$. Then $D$ 
is foliated by smooth, complete, locally strongly convex surfaces asymptotic to $\pa D$, whose Gaussian curvatures with respect to Li's $\alpha$-normal fields (defined using $v$ as the vertical unit vector) are constants ranging from $0$ to $+\infty$. Moreover, the function on $D$ which assigns to each leaf its curvature is $\log$-convex.
\end{manualtheorem}

As a prerequisite for Theorems \ref{thm_2} and \ref{thm_2'}, we establish a classification result for affine hyperspheres with Li-normalization, stated in PDE form and geometric form separately as follows (as before, the parameters are related by $\gamma=\frac{n+2}{\alpha}$):
\begin{thmx}[Simplified version of Theorem \ref{thm_affinesphere1analytic}]\label{thm_3}
For any bounded convex domain $\Omega\subset\R^n$ ($n\geq2$) and any $\gamma>n$, there exists a unique convex generalized solution $w$ to the Dirichlet problem
\begin{equation}\label{eqn_introgamma}
	\begin{cases}
		\det\D^2w=(-w)^{-\gamma}\ \text{ in }\Omega,\\
		w|_{\pa\Omega}=0,
	\end{cases}
\end{equation}
which is in $\C^0(\overline{\Omega})\cap\C^\infty(\Omega)$ and has the gradient blowup property on $\pa\Omega$.
\end{thmx}

Theorem \ref{thm_3} can be seen as a generalization of Theorem \ref{thm_2}, which corresponds to the case $\varphi\equiv 0$. 
Concerning the range of the exponent $\gamma$, we remark that:
\begin{itemize}
\item For $\gamma>1$, the same conclusions in Theorem \ref{thm_2} and Theorem \ref{thm_3} hold if $\Omega$ satisfies both the exterior and the interior sphere conditions at every boundary point, see Theorem \ref{thm_affinesphere1analytic} and Remark \ref{remark_boundary}.
\item No convex solution to Eq.\eqref{eqn_introgamma} with $0<\gamma\leq1$ can have the gradient blowup property throughout $\pa\Omega$, see Theorem \ref{thm_affinesphere1analytic}.
\end{itemize}
On the geometric side, we summarize Theorem \ref{thm_3} and the following remarks in the following statement.

\begin{manualtheorem}{\ref*{thm_3}'}[Theorem \ref{thm_affinesphere1}]\label{thm_4}
	Let $\alpha\in(0,1+\frac{2}{n})$ and $C\subset\R^{n+1}$ be a proper convex cone containing the vertical unit vector $v$. Then for every $c>0$, $C$ contains a unique complete hyperbolic affine hypersphere with $\alpha$-normalization which is centered at the origin $0\in\R^{n+1}$, asymptotic to $\pa C$, and has shape operator $c I$. More generally, for $\alpha\in(0,n+2)$, the same conclusion holds if $C$ has the property that a hyperplane section of it is a convex domain in $\R^n$ satisfying both the exterior and the interior sphere conditions at every boundary point.
On the other hand, if $\alpha\geq n+2$, there does not exist such complete affine hyperspheres.
\end{manualtheorem}
The $\alpha=1$ case is a well known result of Cheng-Yau \cite{chengyau1,chengyau_complete} (with clarifications due to Gigena \cite{gigena}, Li \cite{li_calabi1,li_calabi2} and Sasaki \cite{sasaki}) proving a conjecture of Calabi \cite{calabi}. The generalization here to Li-normalization is first given by Xiong-Yang \cite{xiong-yang} for $\alpha\in(0,1]$ under smoothness and strict convexity assumptions on $\pa\Omega$. While Xiong-Yang's work is based on Lazer-McKenna \cite{lazer-mckenna}, which studies Eq.\eqref{eqn_introgamma} (with a more general right-hand side) under those assumptions, our proof is more similar to \cite[p.67]{chengyau1}, based on barrier functions on simplices and balls.

\subsection*{Organization of the paper} In Sections \ref{sec_li} and \ref{sec_affinehypersphere}, we review the theory of Li-normalization and affine hyperspheres. In Section \ref{sec_cregular}, we review the theory of $C$-regular domains and define \emph{affine $(C,k)$-hypersurfaces}, which form the subclass of constant Gauss-Kronecker curvature hypersurfaces involved in this paper, then we use these notions to give precise statements of the mains results. In Section \ref{sec_legendre}, we explain how Legendre transformation and other basic constructions in convex analysis permit us to reformulate the geometric statements into PDE ones. Finally, we prove the PDE statements in Sections \ref{sec_proof1}, \ref{sec_proof2} and \ref{sec_proof3}.

\subsection*{Acknowledgment} We are grateful to Connor Mooney for helps with the Monge-Amp\`ere theory.

\section{Li-normalization}\label{sec_li}
Let $\A^{n+1}$ denote the $(n+1)$-dimensional real affine space endowed with a translation-invariant volume form $\omega$, and $\R^{n+1}$ denote the underlying vector space. Given a hypersurface $\Sigma\subset\A^{n+1}$ and a transversal vector field $N:\Sigma\to \R^{n+1}$, we recall the following fundamental notions in affine differential geometry (see \eg \cite{MR3382197,MR1311248,MR1200242} for details):

\begin{itemize}
\item The \emph{induced volume form} $\nu$, \emph{induced affine connection} $\nabla$, \emph{affine metric} $h$, \emph{shape operator} $S$ and \emph{transversal connection form} $\tau$ on $\Sigma$ relative to $N$ are defined by the equalities
$$
\nu(X_1,\cdots,X_n)=\omega(X_1,\cdots, X_n,N),
$$
$$
D_XY=\nabla_XY+h(X,Y)N,
$$
$$
D_XN=S(X)+\tau(X)N,
$$
for any tangent vector fields $X,Y$ and $X_1,\cdots,X_n$ on $\Sigma$, where $D$ is the covariant derivative in the ambient affine space.
\item If the affine metric $h$ is  non-degenerate (\ie a pseudo-Riemannian metric), then $\Sigma$ is said to be \emph{non-degenerate} as well, and this property is independent of the choice of $N$. In particular, $h$ is a Riemannian metric if and only if $\Sigma$ is locally strongly convex\footnote{For a $\C^2$ function or hypersurface, by \emph{locally strongly convexity}, we mean positive definiteness of the Hessian.} and $N$ points towards the convex side of $\Sigma$.
\item The transversal vector field $N$ is said to be \emph{equi-affine} if $\tau=0$. 
\item If $\Sigma$ is non-degenerate, $N$ is said to be an \emph{affine normal field} of $\Sigma$ if it is equi-affine and the volume form $\nu_h$ of the metric $h$ coincides with the induced volume form $\nu$. Such an $N$ exists and  is unique up to sign.
\item The \emph{dual covector field} $N^*$ of $N$ is by definition the map $N^*:\Sigma\to\R^{(n+1)*}$ pointwise satisfying
\begin{equation}\label{eqn_dual}
\langle N^*\!,\, N\rangle=1,\quad \langle N^*\!,\, X\rangle=0
\end{equation}
for any tangent vector field $X$ on $\Sigma$, where $\langle\ccdot,\ccdot\rangle$ denotes the pairing between $\R^{(n+1)*}$ and $\R^{n+1}$. In general, $N^*$ cannot determine $N$ because adding a tangent vector field to $N$ does not change its dual. Nevertheless, when $\Sigma$ is non-degenerate, a map $N^*:\Sigma\to\R^{(n+1)*}$ satisfying the second equality in \eqref{eqn_dual} does determine a unique equi-affine transversal vector field $N$ which it is dual to.
\end{itemize}

We henceforth fix a nonzero vector $v\in\R^{n+1}$ as the vertical unit vector. A locally strongly convex hypersurface $\Sigma\subset\A^{n+1}$ is a \emph{locally strongly convex graph} (with respect to $v$) if the constant vector field
$$
Y:\Sigma\to\R^{n+1},\quad Y\equiv v
$$
is transversal to $\Sigma$ and points towards the convex side of $\Sigma$. In this case, we let $\nu_0$ and $h_0$ denote the induced volume form and affine metric on $\Sigma$ relative to $Y$, respectively, and let $\nu_{h_0}$ denote the volume form of the Riemannian metric $h_0$.

\begin{remark}\label{remark_coordinates}
In the literature, one often works with a fixed coordinate system on $\R^{n+1}$ and $\A^{n+1}$, and take $v=(0,\cdots,0,1)$. In this setting, a convex graph is exactly a hypersurface of the form 
$$
\Sigma=\big\{(x,F(x))\mid x\in U\big\},
$$
where $U$ is a domain in $\R^n$ and $F:U\to \R$ is a locally strongly convex function, and one readily checks that
$$
\nu_0=\dx_1\wedge\cdots\wedge\dx_n,\quad h_0=\sum_{1\leq i,j\leq n}\pa^2_{ij}F(x)\dx_i\dx_j,\quad \nu_{h_0}=(\det\D^2F)^\frac{1}{2}\nu_0.
$$
\end{remark}

\begin{propdef}[An-Min Li, see \cite{xu,xu-li-li,xiong-yang,wu-zhao}]\label{propdef}
Let $\Sigma\subset\A^{n+1}$ be a locally strongly convex graph, $N_1:\Sigma\to\R^{n+1}$ be its affine normal field pointing towards the convex side, and $Y:\Sigma\to\R^{n+1}$ be the constant vertical unit vector field as above. Then $N_1$ is exactly the equi-affine transversal vector field on $\Sigma$ with dual (\cf the last bullet point above) given by
$$
N_1^*=\left(\frac{\nu_{h_0}}{\nu_0}\right)^{-\frac{2}{n+2}}Y^*,
$$
where $\frac{\nu_{h_0}}{\nu_0}$ is the density function of $\nu_{h_0}$ with respect to $\nu_0$. More generally, given $\alpha\in\R$, \emph{Li's $\alpha$-normal field}
$N_\alpha:\Sigma\to\R^{n+1}$ is by definition the equi-affine transversal vector field with dual given by
$$
N_\alpha^*=\left(\frac{\nu_{h_0}}{\nu_0}\right)^{-\frac{2\alpha}{n+2}}Y^*.
$$
\end{propdef}
\begin{proof}
Take a coordinate system, let $v$ and $\Sigma$ be as in Remark \ref{remark_coordinates} and parameterize $\Sigma$ by the domain $U\subset\R^n$ through the function $F:U\to\R$. By computations, we may check the following facts:
\begin{itemize}
\item Letting $\D F:=(\pa_1F,\cdots,\pa_nF)$ denote the gradient of $F$, we have
$$
Y^*=(-\D F,1).
$$
\item Any equi-affine transversal vector field 
$N=(N',N^{(n+1)}):U\cong \Sigma\to\R^{n+1}$,
where $N'$ and $N^{(n+1)}$ denote the horizontal and vertical components, respectively, is determined by the component $N^{(n+1)}$ along with the function
$$
\lambda:=N^{(n+1)}-\D F\cdot N'.
$$
In fact, the equi-affine condition implies that $N'=-(\D^2F)^{-1}\,\D\lambda$.
\item
The induced volume form $\nu$ and the affine metric $h$ on $\Sigma$ relative to an equi-affine transversal vector field $N$ as above are related to $\nu_0$ and $h_0$ (\cf Remark \ref{remark_coordinates}) by
\begin{equation}\label{eqn_li1}
\nu=\lambda \nu_0,\ \ h=\lambda^{-1}h_0.
\end{equation}
Moreover, the dual covector field of $N$ is 
\begin{equation}\label{eqn_li2}
N^*=\lambda^{-1}(-\D F,1).
\end{equation}
\end{itemize}
Now take $N=N_1$. By definition of affine normal fields and \eqref{eqn_li1}, we have
$\nu_h=\lambda^{-\frac{n}{2}}\nu_{h_0}=\nu=\lambda \nu_0$, hence
$$
\lambda=\left(\frac{\nu_{h_0}}{\nu_0}\right)^\frac{2}{n+2}.
$$
Combining with \eqref{eqn_li2}, we arrive at the required statement.
\end{proof}

\begin{remark}\label{remark_covariance}
While the affine normal field $N_1$ is covariant with respect to volume-preserving affine transformations of $\A^{n+1}$ and hence independent of the choice of $v$, Li's $\alpha$-normal field $N_\alpha$ with $\alpha\neq1$ does not have this property. This can be seen from the unit hypersphere $\mathbb{S}^n=\{x\in\R^{n+1}\mid |x|=1\}$ or the hyperboloid $\mathbb{H}^n=\big\{x\in\R^{n+1}\,\big|\, x_1^2+\cdots+x_n^2-x_{n+1}^2=-1\big\}$, for which $N_1$ equals the position vector field
$P(x):=x$ up to a sign. When $\alpha\neq1$, one checks that $N_\alpha^*$ (defined with respect to $v=(0,\cdots,0,1)$) is not proportional to $N_1^*$, hence $N_\alpha$ is not proportional to $P$. Due to the affine symmetries of $\mathbb{S}^n$ and $\mathbb{H}^n$, this implies that $N_\alpha$ is not affine-covariant.
\end{remark}

We refer to Lemma \ref{lemma_legendre} below for a more explicit expression of $N_\alpha$ under a particular parametrization of $\Sigma$, as well as expressions of the affine metric and shape operator of $\Sigma$ relative to $N_\alpha$.

When $\alpha=0$, $N_0$ is just the vertical unit vector field $Y$, and we will only consider the $\alpha\neq0$ case below. The following lemma characterizes $N_\alpha$ in a way similar to the definition of affine normal fields:
\begin{lemma}\label{lemma_characterization}
Let $\Sigma$ be as above and suppose $\alpha\neq0$. Then $N_\alpha$ is the unique equi-affine transversal vector field on $\Sigma$ pointing towards the convex side with the following property: if $\nu$ and $h$ are the induced volume form and affine metric on $\Sigma$ relative to $N_\alpha$, respectively, and $\nu_h$ is the volume form of the metric $h$, then the densities of $\nu$ and $\nu_h$ with respect to $\nu_0$ are related by
\begin{equation}\label{eqn_characterization}
\frac{\nu_h}{\nu_0}=\left(\frac{\nu}{\nu_0}\right)^\beta, \ \text{ where }\beta:=\frac{1}{2}\left(\frac{n+2}{\alpha}-n\right).
\end{equation}
\end{lemma}
\begin{proof}
Following the computations in the previous proof, the condition for an equi-affine transversal vector field $N=(N',N^{(n+1)})$ to be Li's $\alpha$-normal field can be written as
\begin{equation}\label{eqn_characterization1}
\lambda=\left(\frac{\nu_{h_0}}{\nu_0}\right)^\frac{2\alpha}{n+2},
\end{equation}
where $\lambda:=N^{(n+1)}-\D F\cdot N'$. On the other hand, by \eqref{eqn_li1}, we have 
$$
\frac{\nu_h}{\nu_0}=\frac{\nu_h}{\nu_{h_0}}\frac{\nu_{h_0}}{\nu_0}=\lambda^{-\frac{n}{2}}\frac{\nu_{h_0}}{\nu_0},\quad \frac{\nu}{\nu_0}=\lambda,
$$
so condition \eqref{eqn_characterization} can be rewritten as
\begin{equation}\label{eqn_characterization2}
\lambda^{-\frac{n}{2}}\frac{\nu_{h_0}}{\nu_0}=\lambda^\beta.
\end{equation}
Conditions \eqref{eqn_characterization1} and \eqref{eqn_characterization2} are clearly equivalent to each other, which proves the required statement.
\end{proof}

\section{Affine hyperspheres and constant curvature with Li-normalization}\label{sec_affinehypersphere}
In this section, we fix a constant $\alpha\neq0$, a vertical unit vector $v\in\R^{n+1}$, let $\Sigma\subset\A^{n+1}$ be a smooth locally strongly convex graph (with respect to $v$), $N_\alpha:\Sigma\to\R^{n+1}$ be Li's $\alpha$-normal field on $\Sigma$, and $S_\alpha$ be the shape operator of $\Sigma$ relative to $N_\alpha$. The following notions are straightforward generalizations of the classical ones:
\begin{itemize}
	\item $\Sigma$ is called a \emph{hyperbolic affine hypersphere with $\alpha$-normalization}\footnote{In the literature \cite{xu,xu-li-li,xiong-yang,wu-zhao} and the introduction, ``$\alpha$-normalization'' has been phrased alternatively as ``Li-normalization''. Here we choose to emphasize the dependence on the parameter $\alpha$ in order to facilitate the statement of some results below.} if 
	$$
	S_\alpha=c I
	$$
	for a constant $c>0$, where $I$ is the identity $(1,1)$-tensor on $\Sigma$. This condition is equivalent to the existence of a point $o\in\A^{n+1}$, the \emph{center} of $\Sigma$, such that $N_\alpha(p)=c\,\overrightarrow{op}$ for all $p\in\Sigma$. The same definition can be made for $c\leq0$ and yields elliptic or parabolic affine hyperspheres, but they are not the concern of this paper.
	\item The \emph{Gauss-Kronecker curvature} (which we simply call \emph{Gaussian curvature} if $n=2$) of $\Sigma$ \emph{with $\alpha$-normalization} is by definition the function
	$$
	\kappa_\alpha:=\det(S_\alpha):\Sigma\to\R.
	$$
\end{itemize}

As a simple fact for any equi-affine transversal vector field $N:\Sigma\to\R^{n+1}$, if the shape operator $S$ of $\Sigma$ relative to $N$ is non-degenerate (\ie $\det(S)\neq0$ everywhere), then $N$ defines an immersion of $\Sigma$ into $\R^{n+1}$.
With this in mind, we can state the following result, which links the above two notions together and generalizes \cite[Prop.\@ 3.5 (2)]{nie-seppi} about the $\alpha=1$ case:
\begin{proposition}\label{prop_link}
Suppose $S_\alpha$ is non-degenerate. Then the following conditions are equivalent:
\begin{enumerate}[label=(\alph*)]
	\item\label{item_link1} The eigenvalues of $S_\alpha$ are all positive, and $\kappa_\alpha$ is a constant.
	\item\label{item_link2} $N_\alpha(\Sigma)\subset\R^{n+1}$ is a hyperbolic affine hypersphere with $\alpha$-normalization, centered at the origin.
\end{enumerate}
Moreover, when these conditions are fulfilled, the affine hypersphere in \ref{item_link2} has shape operator $\kappa_\alpha^{-\frac{\alpha}{n+2}}I$.
\end{proposition}
\begin{proof}
For any equi-affine transversal vector field $N:\Sigma\to\R^{n+1}$, if $\nu$, $h$ and $S$ are the induced volume form, affine metric and shape operator of $\Sigma$ relative to $N$, and $\nu'$, $h'$ and $S'$ are the induced volume form, affine metric and shape operator of $N$ viewed as a centro-affine immersion (\ie relative to the position vector field), then it can be shown that
\begin{equation}\label{eqn_link1}
\nu'=\det(S)\nu,\quad h'(\ccdot,\ccdot)=h(S(\ccdot),\ccdot),\quad S'=I
\end{equation}
(see \cite[Lemma 2.6]{nie-seppi}). A similar computation also shows that if $\nu_0$ and $\nu_0'$ are the induced volume forms of $\Sigma$ and the immersion $N$, respectively, relative to the vertical unit vector field, then
\begin{equation}\label{eqn_link2}
\nu_0'=\det(S)\nu_0.
\end{equation}

Now if condition \ref{item_link1} holds, applying \eqref{eqn_link1} and \eqref{eqn_link2} to $N=N_\alpha$ and using Lemma \ref{lemma_characterization}, we conclude that
\begin{equation}\label{eqn_link3}
\frac{\nu_{h'}}{\nu_0'}=\frac{\det(S)^\frac{1}{2}\nu_h}{\det(S)\nu_0}=\kappa_\alpha^{-\frac{1}{2}}\left(\frac{\nu}{\nu_0}\right)^\beta=\kappa_\alpha^{-\frac{1}{2}}\left(\frac{\nu'}{\nu_0'}\right)^\beta,
\end{equation}
where $\beta:=\frac{1}{2}\left(\frac{n+2}{\alpha}-n\right)$. In general, scaling a transversal vector field by a constant $r$ amounts to scaling the induced volume form by $r$ and the volume form of affine metric by $r^{-\frac{n}{2}}$. Therefore, \eqref{eqn_link3} and Lemma \ref{lemma_characterization} imply that for the hypersurface $N_\alpha(\Sigma)$, scaling the position vector field by 
$$
r:=\kappa_\alpha^{-\frac{1}{2\beta+n}}=\kappa_\alpha^{-\frac{\alpha}{n+2}}
$$ 
yields Li's $\alpha$-normal field. It follows that $N_\alpha(\Sigma)$ is a hyperbolic affine hypersphere with $\alpha$-normalization, centered at the origin, and its shape operator is $r I$. This proves the implication ``\ref{item_link1}$\Rightarrow$\ref{item_link2}'' and the last statement of the proposition.

Conversely, if \ref{item_link2} holds, in view of the expression of $h'$ in \eqref{eqn_link1} and the fact that $h'$ and $h$ are both Riemannian metrics (as $\Sigma$ and $N_\alpha(\Sigma)$ are locally strongly convex; \cf \S \ref{sec_li}), we first conclude that the eigenvalues of $S_\alpha$ are positive. Then, by Lemma \ref{lemma_characterization} and the above scaling argument again, we have
$$
\frac{\nu_{h'}}{\nu_0'}=c\left(\frac{\nu'}{\nu_0'}\right)^\beta
$$
for a constant $c>0$. Using this and equations \eqref{eqn_link1} and \eqref{eqn_link2}, we get
$$
\frac{\nu_h}{\nu_0}=\kappa^\frac{1}{2}_\alpha\frac{\nu_{h'}}{\nu_0'}=c\kappa_\alpha^\frac{1}{2}\left(\frac{\nu'}{\nu_0'}\right)^\beta=c\kappa^\frac{1}{2}_\alpha\left(\frac{\nu}{\nu_0}\right)^\beta.
$$
Therefore, again by Lemma \ref{lemma_characterization}, the assumption that $N_\alpha$ is Li's $\alpha$-normal field implies $c\kappa_\alpha^\frac{1}{2}\equiv1$, hence $\kappa_\alpha$ is a constant. This shows ``\ref{item_link2}$\Rightarrow$\ref{item_link1}'' and completes the proof.
\end{proof}

We henceforth restrict ourselves to \emph{complete}, \ie properly embedded\footnote{This notion is usually referred to as \emph{Euclidean completeness} in the literature on affine differential geometry, in order to make the distinction with \emph{affine completeness}, \ie completeness of the affine metric. We do not study the latter notion in this paper.}, affine hyperspheres with $\alpha$-normalization. Recall that a \emph{convex cone} in $\R^{n+1}$ is by definition a convex domain invariant under positive scalings, and is said to be \emph{proper} if it does not contain any entire affine line. The classification theorem of usual hyperbolic affine hyperspheres, conjectured by Calabi \cite{calabi} and proved by Cheng-Yau \cite{chengyau1,chengyau_complete} (with clarifications due to Gigena \cite{gigena}, Li \cite{li_calabi1,li_calabi2} and Sasaki \cite{sasaki}), roughly claims a $1$-to-$1$ correspondence between proper convex cones in $\R^{n+1}$ and complete affine hyperspheres centered at the origin. Xiong-Yang \cite{xiong-yang} generalized the classification to the setting of $\alpha$-normalization under some extra assumptions. Theorem \ref{thm_4} in the introduction, which we restate below, improves on Xiong-Yang's result:
\begin{theorem}[Theorem \ref{thm_4}]\label{thm_affinesphere1}
	Let $\alpha\in(0,1+\frac{2}{n})$ and $C\subset\R^{n+1}$ be a proper convex cone containing the vertical unit vector $v$. Then for every $c>0$, $C$ contains a unique complete hyperbolic affine hypersphere with $\alpha$-normalization which is centered at the origin $0\in\R^{n+1}$, asymptotic to $\pa C$, and has shape operator $c I$. More generally, for $\alpha\in(0,n+2)$, the same conclusion holds if $C$ has the property that a hyperplane section of it is a convex domain in $\R^n$ satisfying both the exterior and the interior sphere conditions at every boundary point.
On the other hand, if $\alpha\geq n+2$, there does not exist such complete affine hypersphere.
\end{theorem}
Here, a convex domain $\Omega\subset\R^n$ is said to ``satisfy both the exterior and the interior sphere conditions at every boundary point'' if for every $p\in\pa\Omega$ there exist round balls $B_1,B_2\subset\R^n$ such that $B_1\subset\Omega\subset B_2$ and $\pa B_1\cap\pa B_2=\{p\}$.
\begin{remark}\label{remark_scaling}
Using the definitions, it can be shown that if $\Sigma_0$ is an affine hypersphere with $\alpha$-normalization, centered at the origin, whose shape operator is the identity $I$, then the scaled hypersurface $t\Sigma_0$ ($t>0$) is also such an affine hypersphere, with shape operator $t^{-\left(1+\frac{\alpha n}{n+2}\right)}I$. Therefore, the affine hyperspheres given by Theorem \ref{thm_affinesphere1} form a homothetic family.
\end{remark}
We will outline in Section \ref{sec_legendre} the standard process of transforming Theorem \ref{thm_affinesphere1} into an equivalent PDE statement via Legendre transformation, then prove the equivalent statement (\ie Theorem \ref{thm_3} in the introduction) in Section \ref{sec_proof1}.

In view of Proposition \ref{prop_link} and Theorem \ref{thm_affinesphere1}, we now define the following particular class of convex graphs with constant Gauss-Kronecker curvature with $\alpha$-normalization, which we study in the sequel:
\begin{definition}\label{def_ck}
Given $\alpha,k>0$ and a proper convex cone $C\subset\R^{n+1}$ containing the vertical unit vector $v$, a \emph{$(C,k)$-hypersurface with $\alpha$-normalization} is a smooth locally strongly convex graph $\Sigma\subset\A^{n+1}$ such that
\begin{itemize}
	\item the Gauss-Kronecker curvature $\kappa_\alpha$ of $\Sigma$ is constantly $k$;
	\item Li's $\alpha$-normal field $N_\alpha$ of $\Sigma$, viewed as an immersion of $\Sigma$ into $\R^{n+1}$, has image in a complete hyperbolic affine hypersphere with $\alpha$-normalization asymptotic to $\pa C$.
\end{itemize}
\end{definition}
\begin{remark}\label{remark_scaling2}
We only consider the $\alpha<1+\frac{2}{n}$ case of this definition, where the hyperspheres described in the last bullet point uniquely exist up to scaling by Theorem \ref{thm_affinesphere1} and Remark \ref{remark_scaling}. The exactly hypersphere in this scaling family that contains the image of $N_\alpha$ is actually determined by $k$: If we let $\Sigma_C$ denote the one in this family with identity shape operator, then by Proposition \ref{prop_link} and Remark \ref{remark_scaling}, $N_\alpha$ is contained in $k^{\frac{\alpha}{(\alpha+1)n+2}}\Sigma_C$, whose shape operator is $k^{-\frac{\alpha}{n+2}}I$. Therefore, the defining conditions are equivalent to the condition ``$S_\alpha$ is non-degenerate and $N_\alpha$ has image in the scaled affine hypersphere $k^{\frac{\alpha}{(\alpha+1)n+2}}\Sigma_C$''.
\end{remark}

\section{$C$-regular domains and foliation by $(C,k)$-hypersurfaces with $\alpha$-normalization}\label{sec_cregular}
Given a proper convex cone $C\subset\R^{n+1}$, the following definitions introduced in \cite{nie-seppi2,nie-seppi} generalize classical notions from Minkowski geometry:
\begin{itemize}
	\item A subspace $L\subset \R^{n+1}$ of dimension $n$ is said to be \emph{$C$-spacelike} if $L$ meets the closure $\overline{C}$ of $C$ only at the origin $0$, and is said to be \emph{$C$-null} if $L\cap \overline{C}$ is contained in $\pa C$ and is not the single point $0$.
    \item An affine hyperplane in $\A^{n+1}$ is said to be $C$-spacelike/$C$-null if the underlying vector subspace of $\R^{n+1}$ is. Given such a hyperplane $H\subset\A^{n+1}$, we let $\widehat{H}$ denote the closed half-space of $\A^{n+1}$ bounded by $H$ which contains a translation of $C$ (\ie the ``upper half'' of $\A^{n+1}$ cut out by $H$).
    \item A \emph{$C$-regular domain} in $\A^{n+1}$ is by definition an unbounded convex domain $D\subset \A^{n+1}$ of the form 
    $$
    D=\interior\bigcap_{H\in\mathcal{S}}\widehat{H}
    $$ 
    where $\interior$ denotes the interior and $\mathcal{S}$ is a collection of $C$-null hyperplanes. If $\mathcal{S}$ consists of all $C$-null hyperplanes $H$ such that $\widehat{H}$ contains some set $E\subset\A^{n+1}$, then $D$ is called the $C$-regular domain \emph{generated} by $E$. 
    \item A \emph{$C$-convex hypersurface} is by definition an open subset $\Sigma$ of the boundary of some convex domain $U\subset\A^{n+1}$ such that any supporting hyperplane $H$ of $U$ at any point of $\Sigma$ is $C$-spacelike with $U\subset\widehat{H}$. 
    So $\Sigma$ is complete if it is the entire $\pa U$. The simplest examples are just $C$-spacelike hyperplanes. The main class of examples that we study are affine $(C,k)$-hypersurfaces from Definition \ref{def_ck}, which can be shown to be $C$-convex.
\end{itemize}  

For the future light cone $C_0$ in the Minkowski space $\R^{n,1}$, $C_0$-regular domains are classically known as \emph{regular domains} or \emph{domains of dependence}, and play an important role in the study of \emph{globally hyperbolic flat spacetimes} from mathematical relativity, where one considers a discontinuous isometric action on such a domain with quotient a Lorentzian manifold of the form $M\times\R$, where $M$ is a compact $n$-manifold and $\R$ is the time direction (see \eg \cite{MR2110829,bbz,MR2170277,mess}). Mess \cite{mess} discovered a profound link between the $n=2$ case and Teichm\"uller Theory. Some of Mess' results are generalized from $C_0$ to general $C$ in \cite{nie-seppi2}. 

The main results of this paper can be viewed as generalizations of Theorem \ref{thm_affinesphere1} above, with the cone $C$ replaced by $C$-regular domains, and the affine hyperspheres replaced by $(C,k)$-hypersurfaces. We first give a statement which holds for any $n\geq2$:
\begin{theorem}[Extended version of Theorem \ref{thm_2'}]\label{thm_ck1}
Let $0<\alpha<1+\frac{2}{n}$, $C\subset\R^{n+1}$ be a proper convex cone  containing the vertical unit vector $v$, and $D\subset\A^{n+1}$ be a $C$-regular domain satisfying the following condition: for every $C$-null subspace $L\subset\R^{n+1}$, there exists a translation $C'\subset\A^{n+1}$ of the cone $C$ such that $C'$ contains $D$ and the supporting hyperplane $L'$ of $C'$ parallel to $L$ is an asymptotic hyperplane\footnote{An \emph{asymptotic hyperplane} of an unbounded convex domain $U\subset\A^{n+1}$ is by definition a hyperplane $H$ disjoint from $U$ such that any hyperplane $H'$ obtained by moving $H$ parallelly towards $U$ intersects $U$ along an unbounded set.}
 of $D$. Then $D$ contains, for every $k>0$, a unique complete $(C,k)$-hypersurface with $\alpha$-normalization $\Sigma_k$ which generates $D$. Moreover, the following statements hold: 
\begin{itemize} 
 \item For each $k$, the distance from $x\in\Sigma_k$ to $\pa D$ tends to $0$ as $x$ tends to infinity in $\Sigma_k$.
 \item $(\Sigma_k)_{k>0}$ is a foliation of $D$ (\ie the $\Sigma_k$'s are disjoint and their union is $D$).
 \item
 The function $K:D\to \R$ given by $K|_{\Sigma_k}=\log k$ is convex. 
\end{itemize}
\end{theorem}
\begin{remark}\label{remark_ck1}
\begin{enumerate}
\item The assumption on $D$ is not satisfied, for example, when $D=C\cap \interior\,\widehat{L_0'}$ for some $C$-null subspace $L_0$ and a $C$-null hyperplane $L_0'$ obtained by moving $L_0$ towards $C$. 
\item	
The condition that ``the distance from $x\in\Sigma$ to $\pa D$ tends to $0$ as $x$ goes to infinity in $\Sigma$'' is phrased as ``$\Sigma$ is asymptotic to $\pa \Sigma$'' in \cite[\S 5.3]{nie-seppi}, and it is explained therein that for a general complete $C$-convex hypersurface $\Sigma$, this condition is strictly stronger then ``$\Sigma$ generates $D$''.
\end{enumerate}
\end{remark}
Like Theorem \ref{thm_affinesphere1}, we will transform Theorem \ref{thm_ck1} and the other two results below into equivalent PDE statements via the construction in the next section, and then give proofs in Sections \ref{sec_proof2} and \ref{sec_proof3}.
The reader may check the PDE version Theorem \ref{thm_ck1analytic} of Theorem \ref{thm_ck1} (see also Remark \ref{remark_boundary}) for a more transparent interpretation of the assumption on $D$. All these results improve on the work of Wu-Zhao \cite{wu-zhao}.

When $n=2$, under some mild restrictions on $\alpha$ and $C$, we may extend Theorem \ref{thm_ck1} to any $C$-regular domain $D$ which is proper (\ie not containing any entire affine line), without any extra assumption:
\begin{theorem}[Extended version of Theorem \ref{thm_1'}]\label{thm_ck2}
Suppose $0<\alpha\leq 1$ and let $C\subset\R^3$ be a proper convex cone containing the vertical unit vector such that a plane section $\Omega\subset\R^2$ of the dual cone $C^*\subset\R^{3*}$ satisfies the exterior circle condition. Then the conclusions of Theorem \ref{thm_ck1} hold for any proper $C$-regular domain $D\subset\A^3$.  
\end{theorem}
This is a generalization of our recent work \cite{nie-seppi} about the $\alpha=1$ case.

In Prop.\@ E and Cor.\@ 8.6 of \cite{nie-seppi}, we showed in the $\alpha=1$ case that if $\Omega$ is not strictly convex, then the conclusion does not hold for certain examples of $D$. We believe there also exist such examples where $\Omega$ is strictly convex but does not satisfy the exterior circle condition. On the other hand, the necessity of the condition $\alpha\leq1$ is shown by the following result:
\begin{proposition}\label{prop_alpha1}
Suppose $n=2$ and $\alpha>1$. Let $C_0$ be the future light cone in the Minkowski space $\R^{2,1}$ and $D$ be a triangular cone bounded by three null planes, which is a $C_0$-regular domain. Then $D$ does not contain any complete $(C_0,k)$-surface with $\alpha$-normalization which generates $D$.
\end{proposition}

\section{Legendre transformation}\label{sec_legendre}
We henceforth identify $\A^{n+1}$ with $\R^{n+1}$ by choosing a base point, and fix a coordinate system on $\R^{n+1}$ under which the vertical unit vector is $v=(0,\cdots,0,1)$. We usually write a point of $\R^{n+1}$ as $X=(x,\xi)$, where $x\in\R^n$ and $\xi\in\R$ are the horizontal and vertical components, respectively.

Given a proper convex cone $C\subset\R^{n+1}$ containing $v$, an important bounded convex domain $\Omega\subset\R^{n}$ associated with $C$ is the section of the dual cone $C^*:=\big\{X\in \R^{n+1}\,\big|\, X\cdot Y<0,\,\forall Y\in\overline{C}\setminus\{0\}\big\}$\footnote{Conceptually, $C^*$ is a convex cone in the dual vector space $\R^{(n+1)*}$, but here we use the standard inner product ``$\ \cdot\ $'' to identify $\R^{(n+1)*}$ with $\R^{n+1}$. Also, the definition of $C^*$ in the literature (\eg \cite{nie-seppi}) sometimes differs from the one here by a sign.} by the hyperplane $\R^{n}\times\{-1\}$. That is, $\Omega$ is defined by
$$
C^*=\big\{t(x,-1)\,\big|\,x\in\Omega,\,t>0\big\}.
$$
As an alternative definition, if we suppose $C=\big\{t(y,1)\,\big|\,y\in C_1,\,t>0\big\}$, namely $C_1$ is the section of $C$ by the hyperplane $\R^n\times\{1\}$, then we have
$$
\Omega=\big\{x\in\R^n\,\big|\,x\cdot y<1,\,\forall y\in\overline{C}_1\big\}.
$$
That is, $\Omega$ is the dual convex domain of $C_1$ in the sense of Sasaki \cite{sasaki}.

The significance of $\Omega$ lies in the fact that all the geometric objects in $\A^{n+1}$ introduced above with respect to $C$, such as $C$-null/$C$-spacelike hyperplanes, $C$-regular domains, $(C,k)$-hypersurfaces, \emph{etc.}, can be interpreted in terms of convex functions on $\Omega$ through Legendre transformation, which we now explain.

Let us first recall some facts about convex function. Let $\LC(\R^n)$ denote the space of $\R\cup\{+\infty\}$-valued, lower semicontinuous, convex functions on $\R^n$ that are not constantly $+\infty$. Given $u\in\LC(\R^n)$, if the effective domain $\dom{u}:=\{x\mid u(x)<+\infty\}$ has nonempty interior $U:=\interior\dom{u}$, then the values of $u$ on $\pa U$ (hence the values on the whole $\R^n$) are determined by the restriction $u|_U$, because we have
\begin{equation}\label{eqn_boundary value}
u(p)=\liminf_{U\ni x\to p}\,u(x)=\lim_{s\to0^+}u((1-s)p+sx)\in(-\infty,+\infty]
\end{equation}
for any $p\in\pa U$ and $x\in U$ (see \cite[\S 4.1]{nie-seppi}). 

Therefore, given a convex domain $U\subset \R^n$ and a convex function $u:U\to\R$, we define the \emph{boundary value} $u|_{\pa U}$ of $u$ as the function on $\pa U$ whose value at any $p\in\pa U$ is the liminf or the limit in \eqref{eqn_boundary value}. The extension of $u$ to $\R^n$ given by $u|_{\pa U}$ and by setting $u=+\infty$ outside of $\overline{U}$ is an element of $\LC(\R^n)$. This gives a canonical way of viewing every convex function on a convex domain as an element of $\LC(\R^n)$.

Given $p\in\pa U$, we say that $u$ has \emph{infinite slope} (or \emph{infinite inner derivatives}) at $p$ if either $u(p)=+\infty$ or $u(p)$ is finite but
\begin{equation}\label{eqn_inner}
\lim_{s\to0^+}\frac{u((1-s)p+sx)-u(p)}{s}=-\infty
\end{equation}
for some $x\in U$ (by convexity of $u$, the fraction decreases as $s$ decreases to $0$, hence the limit exists in $[-\infty,+\infty)$). 
We refer to \cite[\S 4]{nie-seppi} for the following fundamental facts about this notion:
\begin{itemize}
\item It is independent of the choice of $x\in U$. That is, if \eqref{eqn_inner} holds for one $x$, then it holds for all.
\item The following conditions are equivalent:
\begin{itemize}
	\item $u$ has infinite slope at $p\in\pa U$;
	\item $u$ does not have any \emph{subgradient} at $p$ (see \eg \cite[\S 4.4]{nie-seppi} for the definition), or in other words, the graph of $u$ does not have any non-vertical supporting hyperplane at $p$;
	\item for any sequence $(x_i)_{i=1,2,\cdots}$ in $U$ tending to $p$ such that $u$ is differentiable at every $x_i$, we have $\|\D u(x_i)\|\to+\infty$.
\end{itemize}
\end{itemize}
In particular, when $u\in\C^1(U)$, the condition that $u$ has infinite slope at every point of $\pa U$ is equivalent to the gradient blowup condition ``$\|\D u(x)\|\to+\infty$ as $x\in U$ tends to $\pa U$''.


We will consider some particular classes of functions $u\in\LC(\R^n)$ with $\dom{u}$ contained in the closure $\overline{\Omega}$ of the above convex domain $\Omega$. The first consists of \emph{convex envelopes} of functions on $\pa\Omega$, defined as follows: given a function $\psi:\pa\Omega\to\R\cup\{+\infty\}$ which is bounded from below and not constantly $+\infty$, the convex envelope $\env{\psi}\in\LC(\R^n)$ of $\psi$ is 
$$
\env{\psi}(x):=\sup\big\{a(x)\,\big|\, \text{ $a:\R^n\to\R$ is an affine function with $a|_{\pa\Omega}\leq \psi$}\big\}.
$$
As shown in \cite[\S 4.3]{nie-seppi}, the assignment $\phi\mapsto\env{\phi}$ is a bijection from the set of functions
$$
\LC(\pa\Omega):=
\left\{u:\pa\Omega\to\R\cup\{+\infty\}\ \Bigg|\ \parbox[l]{6.8cm}{$u$ is lower semicontinuous and restricts to a convex function on any line segment in $\pa\Omega$.}\right\}~
$$ 
to the subset of $\LC(\R^n)$ consisting of all convex envelopes, with inverse the restriction map $u\mapsto u|_{\pa\Omega}$. We will make use of the following properties of the convex envelope $\env{\phi}$ of $\phi\in\LC(\pa\Omega)$:
\begin{itemize}
	\item $\dom{\env{\phi}}$ is the convex hull of $\dom{\phi}\subset\pa\Omega$ (see \cite[Prop.\@ 4.8]{nie-seppi}).
	\item For any $u\in\LC(\R^n)$ with $u|_{\pa\Omega}\leq\phi$, we have $u\leq\env{\phi}$ throughout $\R^n$ (see \cite[Cor.\@ 4.5]{nie-seppi}).
	\item Suppose $u:\Omega\to\R$ is a convex function with $u|_{\pa\Omega}\leq\phi$. Then either of the following conditions is sufficient for the strict inequality $u<\env{\phi}$ to hold in $\Omega$:
	\begin{itemize}
		\item $u$ is strictly convex;
		\item $u$ has infinite slope at every point of $\pa\Omega$.
	\end{itemize}
This can be proved by using \cite[Lemma 4.9]{nie-seppi}.
    \item If $\phi$ is $\R$-valued, so that $\env{\phi}$ is also $\R$-valued in $\Omega$, then the Monge-Amp\`ere measure of $\env{\phi}$ is the zero measure (see \cite[Thm.\@ 1.5.2]{gutierrez}).
\end{itemize}

Two other classes of $u\in\LC(\R^n)$ that we will consider are
$$
\SC(\Omega):=\left\{ u\in\LC(\R^n)\ \big|\ 
\text{$\dom{u}\subset\overline{\Omega}$, $u$ does not admit any subgradient at any point of $\pa\Omega$}
\right\}
$$
and the subset of $\SC(\Omega)$ given by
$$
\SC_0(\Omega):=\left\{u\in\LC(\R^n)\ \Bigg|\
\parbox[l]{8cm}{the interior $U$ of $\dom{u}$ is nonempty and contained in $\Omega$; $u$ is smooth and locally strongly convex in $U$, and has infinite slope at every point of $\pa U$}
\right\}.
$$

The \emph{Legendre transform} of $u\in\LC(\R^n)$ is by definition the function $u^*\in\LC(\R^n)$ given by
$$
u^*(y):=\sup_{x\in\R^n}(x\cdot y-u(x)).
$$
It is a fundamental fact that the Legendre transformation $u\mapsto u^*$ is an involution on $\LC(\R^n)$ (see \eg \cite[\S 4.5]{nie-seppi}). 
If $u$ is an $\R$-valued convex function only defined on a convex domain $U\subset\R^n$, we can still define $u^*$, either by changing the range of $x$ in the supremum into $U$, or by viewing $u$ as an element of $\LC(\R^n)$ through the canonical extension mentioned above. We will need the following properties of $u^*$: 
\begin{itemize}
	\item $u_1\leq u_2$ on $\R^n$ if and only if $u_1^*\geq u_2^*$ on $\R^n$;
	\item if $u$ is differentiable at a point $x$, then the value of $u^*$ at the gradient $\D u(x)$ is given by
	$$
	u^*(\D u(x))=x\cdot \D u(x)-u(x).
	$$
\end{itemize}
Concerning the last property, we further note that if $u$ is strictly convex and $\C^1$ in a convex domain $U\subset\R^n$, then the gradient map $x\mapsto \D u(x)$ is a homeomorphism from $U$ to the image $\D u(U)$. The property implies that the graph of $u^*$ over $\D u(U)$ is $\big\{\big(\D u(x),\, x\cdot\D u(x)-u(x)\big)\in\R^{n+1}\,\big|\,x\in U \big\}$. Also note that $\D u(U)$ is the whole $\R^n$ if and only if $U=\interior\dom{u}$ and $u$ has infinite slope at every point of $\pa U$.  

With the above notions and facts in mind, we can now formulate the relationship between convex functions on $\Omega$ and the geometric objects in $\A^{n+1}$ as the following theorem. Here and below, given a function $f:\R^n\supset E\to\R\cup\{+\infty\}$, we let 
$$
\gra{f}:=\{(x,\xi)\in\R^{n+1}\mid x\in E,\, \xi=f(x)\},\quad \sepi{f}:=\{(x,\xi)\in\R^{n+1}\mid x\in E,\,\xi>f(x)\}
$$ 
denote the graph and strict epigraph of $f$, respectively.
\begin{theorem}\label{thm_dictionary}
Suppose $u\in\LC(\R^n)$ and $\dom{u}\subset\overline{\Omega}$. Then for each row of Table \ref{table}, $u$ has the given property on the left if and only if the graph of the Legendre transform $u^*$ fulfills the description on the right.
\begin{table}[ht]
	\begin{TAB}[4pt]{|c|c|c|}{|c|c|c|c|c|c|c|c|}
&		\textbf{property of $u$}& \textbf{nature of the graph of $u^*$}\\
\circled{1}&		$u=+\infty$ except at a single point of $\Omega$ (resp.\@ $\pa\Omega$) & a $C$-null (resp. $C$-spacelike) hyperplane\\
\circled{2}&		$u$ is an affine function on $\overline{\Omega}$ & the boundary of a translation of $C$\\
\circled{3}&		$u=\env{\phi}$ for some $\phi\in\LC(\pa\Omega)$  & 
\parbox[l]{5.8cm}{the boundary of a $C$-regular domain $D$ (\ie $D=\sepi{u^*}$ is a $C$-regular domain)}
\\
\circled{4}&$u\in\SC(\Omega)$& complete $C$-convex hypersurface\\
\circled{5}&	$u\in\SC_0(\Omega)$ & \parbox[l]{5.1cm}{complete, smooth, locally strongly convex, $C$-convex hypersurface} \\
\circled{6}&	$u\in\SC_0(\Omega)$, $u|_{\pa\Omega}=0$, $\det\D^2u=(-u)^{-\frac{n+2}{\alpha}}$ in $\Omega$& \parbox[l]{5.1cm}{complete hyperbolic affine hypersphere with $\alpha$-normalization, centered at $0$ and asymptotic to $\pa C$}\\
\circled{7}&	    
\parbox[l]{5.8cm}{$u\in\SC_0(\Omega)$, $\det\D^2u=k^{-\frac{n+2}{(\alpha+1)n+2}}(-w)^{-\frac{n+2}{\alpha}}$ in the interior of $\dom{u}$ for some $w\in\SC_0(\Omega)$ satisfying the condition in \circled{6}}&
\parbox[l]{4.2cm}{complete $(C,k)$-hypersurface with $\alpha$-normalization
} 
	\end{TAB}
	\caption{Correspondence between convex functions on $\Omega$ and objects in $\A^{n+1}$}
	\label{table}
\end{table}
Moreover, the following statements about the rows \circled{3} and \circled{4} hold true:
\begin{enumerate}[label=(\arabic*)]
	\item\label{item_dictionary0} Given $\phi\in\LC(\pa\Omega)$, the $C$-regular domain $D=\sepi{\env{\phi}^*}$ is proper if and only if the convex set $\dom{\env{\phi}}$ (which is the convex hull of $\dom{\phi}$) has nonempty interior.
	\item\label{item_dictionary1} Suppose $u\in\SC(\Omega)$ and denote $\phi:=u|_{\pa\Omega}\in\LC(\pa\Omega)$. Then the hypersurface $\Sigma:=\gra{u^*}$ generates the domain $D:=\sepi{\env{\phi}^*}$ and is contained in $D$. In this case, the distance from $x\in\Sigma$ to $\pa D$ tends to $0$ as $x$ goes to infinity in $\Sigma$ (\cf Remark \ref{remark_ck1}) if and only if the following conditions are satisfied:
	\begin{itemize}
    \item
	$\dom{u}$ coincides with $\dom{\env{\phi}}$;
	\item
    $\env{\phi}(x)-u(x)\to0$ as $x\in U:=\interior\dom{u}$ tends to $\pa U$. 
    \end{itemize}
    \item\label{item_dictionary2} Let $(u_t)_{t\in\R}$ be a one-parameter family in
     $\SC(\Omega)$ with the same boundary value $\phi=u_t|_{\pa\Omega}\in\LC(\pa\Omega)$, satisfying the following conditions:
     \begin{itemize}
     	\item for any $x\in\R^n$ and $t_1<t_2$, we have $u_{t_1}(x)\leq u_{t_2}(x)$, and the inequality is strict if $u_{t_1}$ admits a subgradient at $x$;
     	\item as $t\to+\infty$, $u_t$ pointwise converges to $\env{\phi}$;
     	\item $u_t$ is not bounded from below uniformly in $t$;
     	\item for every fix $x\in\R^n$, the function $\R\to\R\cup\{+\infty\}$, $t\to u_t(x)$ is concave.
     \end{itemize}
 Then the hypersurfaces $\Sigma_t:=\gra{u^*_t}$, $t\in\R$ form a foliation of the domain $D:=\sepi{\env{\phi}^*}$ and has the property that the function $K:D\to \R$ given by $K|_{\Sigma_t}= t$ is convex. Meanwhile, the above conditions are also necessary for the $\Sigma_t$'s to form a foliation of $D$ with that property.
\end{enumerate}
\end{theorem}
\begin{remark}
A special instance of row \circled{3} is the whole affine space $\A^{n+1}$ as a $C$-regular domain, whose corresponding $u$ is the constant function $+\infty$ (although $u$ is not in $\LC(\R^n)$, the Legendre transform $u^*$ is well defined and equals $-\infty$). Consequently, as the simplest example of statement \ref{item_dictionary2}, given a point $x_0\in\Omega$, the family $(u_t)_{t\in\R}$ given by $u_t(x_0)=t$ and $u_t=+\infty$ on $\R^n\setminus\{x_0\}$ corresponds to the foliation of $\A^{n+1}$ by the translations of a $C$-spacelike hyperplane.
\end{remark}

The equivalences in the rows \circled{1}$\sim$\circled{5} and statements  \ref{item_dictionary0}, \ref{item_dictionary1} and \ref{item_dictionary2} of the theorem are essentially the content of \cite[\S 5]{nie-seppi}, so we omit the proof of these parts. The affine differential geometry computations leading to the rows \circled{6} and \circled{7} have been done in \cite{xiong-yang, wu-zhao} and can be summarized as follows:
\begin{proposition}\label{prop_legendre}
	Let $U\subset\R^n$ be a convex domain and $u\in\C^\infty(U)$ be a locally strongly convex function. Let $\Sigma$ denote the graph of $u^*$ over $\D u(U)$, \ie
	$$
	\Sigma:=\gra{u^*|_{\D u(U)}}=\big\{\big(\D u(x),\,x\cdot\D u(x)-u(x)\big)\,\big|\,x\in U \big\}\subset\R^{n+1}.
	$$
Then the following statements hold:
	\begin{enumerate}[label=(\arabic*)]
		\item\label{item_legendre1} $\Sigma$ is a hyperbolic affine hypersphere with $\alpha$-normalization, centered at the origin $0\in\R^{n+1}$, and has shape operator $cI$, if and only if $u$ satisfies
		$$
		\det\D^2u=(-c\,u)^{-\frac{n+2}{\alpha}}.
		$$
		\item\label{item_legendre2} Let $N_\alpha:\Sigma\to\R^{n+1}$ be Li's $\alpha$-normal field of $\Sigma$,  $\widetilde{u}\in\C^\infty(U)$ be another locally strongly convex function 
        and $\widetilde{\Sigma}$ be the graph of $\widetilde{u}^*$ over $\D\widetilde{u}(U)$. Then $N_\alpha$ is an immersion with  image contained in $\widetilde{\Sigma}$ if and only if $u$ satisfies
		$$
		\det\D^2u=(-\widetilde{u})^{-\frac{n+2}{\alpha}}.
		$$
		Moreover, in this case, $N_\alpha$ is actually a diffeomorphism from $\Sigma$ to $\widetilde{\Sigma}$.
	\end{enumerate}
\end{proposition}
The equivalence between the properties of $u$ and $\gra{u^*}$ in row \circled{6} of Theorem \ref{thm_dictionary} follows immediately from part \ref{item_legendre1} of the proposition and row \circled{5}, whereas the the equivalence in row \circled{7} follows from part \ref{item_legendre2} in combination with Remark \ref{remark_scaling2}.

For the sake of exposition, we give a proof of Prop.\@ \ref{prop_legendre} using the lemma below, which provides an explicit expressions of Li's $\alpha$-normal field $N_\alpha$, as well as the affine metric and shape operator relative to $N_\alpha$, although the latter two are not needed in Prop.\@ \ref{prop_legendre}. Here, for any $u\in\C^1(U)$ (not necessarily convex) we call the map 
$U\to\R^{n+1}$, $x\mapsto \big(\D u(x),\,x\cdot\D u(x)-u(x)\big)
$ the \emph{Legendre map} of $u$.
\begin{lemma}\label{lemma_legendre}
	Let $U$ and $u$ be as in Proposition \ref{prop_legendre} and $f:U\to\R^{n+1}$ be the Legendre map of $u$, viewed as a convex hypersurface embedding. Put
	$w_\alpha:=-(\det\D^2 u)^{-\frac{\alpha}{n+2}} \in\mathsf{C}^\infty(U)$.
	Then the Legendre map 
	$$
	N_\alpha(x):=\big(\D w_\alpha(x),\,x\cdot \D w_\alpha(x)-w_\alpha(x)\big)
	$$ 
	of $w_\alpha$ is Li's $\alpha$-normal field for the embedding $f$ (with respect to the vertical unit vector $v=(0,\cdots,0,1)$). Moreover, the affine metric $h_\alpha$ and the shape operator $S_\alpha$ of $f$ relative to $N_\alpha$ have the following matrix expressions (with respect to the frame $(\pa_1,\cdots,\pa_n)$):
	$$
	h_\alpha=-\frac{1}{w}\D^2u,\quad S_\alpha=(\D^2 u)^{-1}\D^2w_\alpha.
	$$
\end{lemma}
\begin{proof}
By standard computations, one finds that for a general vector field $N=(N',N^{(n+1)}):U\to\R^{n+1}$ transversal to $f$ (where $N':U\to\R^n$ and $N^{(n+1)}:U\to\R$ are respectively the horizontal and vertical components), if we define
$$
w:U\to\R,\quad w(x):=x\cdot N'(x)-N^{(n+1)}(x),
$$
then the induced volume form $\nu$, affine metric $h$, shape operator $S$ and transversal connection form $\tau$ of $f$ relative to $N$ are given by
$$
\nu=-w\det\D^2u\dx_1\wedge\cdots\wedge\dx_n,\quad \tau=-\frac{1}{w}\left(\dif N^{(n+1)}-x\cdot\dif N'\right),
$$
$$
h=-\frac{1}{w}\D^2u,\quad S=(\D^2u)^{-1}\left(\pa_1N',\cdots,\pa_nN'\right)-N'\,(\tau_1,\cdots,\tau_n),
$$
where in the last matrix expression, $N'$ is understood as a row vector of functions, and $\tau_1,\cdots,\tau_n$ are the functions such that $\tau=\tau_1\dx_1+\cdots+\tau_n\dx_n$.

The expression of $\tau$ implies that if $N$ is equi-affine (\ie $\tau=0$), then we have
$$
\D w(x)=\D\left(x\cdot N'(x)-N^{(n+1)}(x)\right)=N'(x),
$$
which implies that $N$ is the Legendre map of $w$. Also, in this case, we have
$$
S=(\D^2u)^{-1}\left(\pa_1N',\cdots,\pa_nN'\right)=(\D^2u)^{-1}\D^2w.
$$
Thus, in order to prove the required statements, it now remains to be shown that if $w=-(\det\D^2u)^{-\frac{\alpha}{n+2}}$ then $N$ is the $\alpha$-normal field of $f$.
 
To this end, we use the characterization of the $\alpha$-normal field in Lemma \ref{lemma_characterization}. By a simple computation, the induced volume form $\nu_0$ relative to the constant vertical unit vector field is given by
$$
\nu_0=\det\D^2u\dx_1\wedge\cdots\wedge\dx_n.
$$ 
So the density functions of $\nu_h$ and $\nu$ with respect to $\nu_0$ are given by
	$$
	\frac{\nu_h}{\nu_0}=\frac{(-w)^{-\frac{n}{2}}(\det\D^2 u)^\frac{1}{2}}{\det\D^2 u}=(-w)^{-\frac{n}{2}}(\det\D^2 u)^{-\frac{1}{2}},\qquad
	\frac{\nu}{\nu_0}=\frac{-w\det\D^2 u}{\det\D^2 u}=-w.
	$$
By Lemma \ref{lemma_characterization}, the condition for $N$ to be the $\alpha$-normal field is
$$
(-w)^{-\frac{n}{2}}(\det\D^2 u)^{-\frac{1}{2}}=(-w)^\beta, \ \text{ where }\beta:=\frac{1}{2}\left(\frac{n+2}{\alpha}-n\right).
$$
This is equivalent to the required equality $w=-(\det\D^2u)^{-\frac{\alpha}{n+2}}$, so the proof is completed.
\end{proof}

For the proof of Prop.\@ \ref{prop_legendre}, we will also need the following facts about the Legendre maps $f_1$ and $f_2$ of two functions $u_1,u_2\in\C^1(U)$ (not necessarily convex):
	\begin{itemize}
		\item We have $f_1=cf_2$ for a constant $c\neq0$ if and only if $u_1=cu_2$.
		\item Suppose $u_1$ is strictly convex and $f_2$ is an immersion of $U$ into $\R^{n+1}$ with image contained in the image of $f_1$. Then $u_1=u_2$.
	\end{itemize}
The first can be checked directly using the definition. For the second statement, first note that the assumption $f_2(U)\subset f_1(U)$ means every supporting hyperplane of the graph $\gra{u_2}$ is also a supporting hyperplane of $\gra{u_1}$ (for example, $u_2$ can be the affine function whose graph is a supporting hyperplane of $\gra{u_1}$). On the other hand, it can be shown that $f_2$ is an immersion with locally strictly convex image if and only if $u_2$ is a locally strictly convex function. Thus, we can only have $u_1=u_2$.

\begin{proof}[Proof of Proposition \ref{prop_legendre}]
\ref{item_legendre1} Let $f$ be the Legendre map of $u$, so that $\Sigma$ is the image of $f$. By definitions, $\Sigma$ is as described in the statement if and only if Li's $\alpha$-normal field $N_\alpha:U\to \R^{n+1}$ of $f$ satisfies 
	\begin{equation}\label{eqn_prooflegendre1}
	N_\alpha=cf.
	\end{equation}
	Lemma \ref{lemma_legendre} says that $N_\alpha$ is the Legendre map of $w_\alpha:=-(\det\D^2u)^{-\frac{\alpha}{n+2}}$. Therefore, by the first statement before the proof, Eq.\eqref{eqn_prooflegendre1} is equivalent to $w_\alpha=cu$, which is in turn equivalent to the required equation.
	
	\ref{item_legendre2} Similarly as above, the condition ``$N_\alpha$ is an immersion with  image contained in $\widetilde{\Sigma}$'' means that the Legendre map of $w_\alpha$ is an immersion with image contained in the image of the Legendre map of $\widetilde{u}$. By the second statement before the proof, this is equivalent to $w_\alpha=\widetilde{u}$, which is in turn equivalent to the required equation and implies the ``Moreover'' statement.
\end{proof}

\section{Monge-Amp\`ere problem for affine hyperspheres}\label{sec_proof1}
In the rest of the paper, we translate the geometric statements above, namely Theorems \ref{thm_affinesphere1}, \ref{thm_ck1}, \ref{thm_ck2} and Proposition \ref{prop_alpha1}, into equivalent PDE statements via Theorem \ref{thm_dictionary}, and then give proofs. We will make use of basic notions in Monge-Amp\`ere theory, such as \emph{Monge-Amp\`ere measure} and \emph{convex generalized solutions} (see \eg \cite{figalli,gutierrez,trudwang}), as well as the notions and facts about convex functions reviewed in Section \ref{sec_legendre}. In particular, recall from Section \ref{sec_legendre} that the gradient blowup property ``$\|\D u(x)\|\to+\infty$ as $x\in\Omega$ tends to $\pa\Omega$'' for a convex function $u:\Omega\to\R$ is equivalent to the \emph{infinite slope property}, which we use in all the statements below.

In this section, we treat Theorem \ref{thm_affinesphere1}, for which it suffices to consider the $c=1$ case by Remark \ref{remark_scaling}. By Theorem \ref{thm_dictionary}, this case is equivalent to the following result, which is stated in the introduction as Theorem \ref{thm_3} and the subsequent discussion:
\begin{theorem}[Extended version of Theorem \ref{thm_3}, PDE version of Theorem \ref{thm_affinesphere1}]\label{thm_affinesphere1analytic}
Let $\Omega$ be a bounded convex domain in $\R^n$ ($n\geq2$) and consider the Dirichlet problem
\begin{equation}\label{eqn_affinesphere2}
	\begin{cases}
	\det\D^2w=(-w)^{-\gamma}\ \text{ in }\Omega,\\
	w|_{\pa\Omega}=0.
	\end{cases}
\end{equation}
\begin{enumerate}[label=(\arabic*)]
\item\label{item_analytic1}
If $\gamma>n$, then Eq.\eqref{eqn_affinesphere2} has a unique convex generalized solution, which is smooth in $\Omega$ and has infinite slope at every point of $\pa\Omega$.
\item\label{item_analytic2}
If $\gamma>1$ and $\Omega\subset\R^n$ satisfies both the exterior and the interior sphere conditions at every boundary point, then the same conclusions as \ref{item_analytic1} holds.
\item\label{item_analytic3}
If $0<\gamma\leq1$, the Eq.\eqref{eqn_affinesphere2} does not have any convex generalized solution with infinite slope at every boundary point.
\end{enumerate}
\end{theorem}
\begin{remark}
	The proof of Part \ref{item_analytic3} shows more specifically that if $\Omega$ satisfies the exterior sphere condition at some $p\in\pa\Omega$, then any convex function $w$ satisfying
	$\det\D^2w\leq (-w)^{-\gamma}$ in the generalized sense with $0<\gamma\leq1$ must have finite slope at $p$.
\end{remark}

The exponent $-\gamma$ in Eq.\eqref{eqn_affinesphere2}  corresponds to the exponent $-\frac{n+2}{\alpha}$ in the equations displayed in Theorem \ref{thm_dictionary}, so the conditions $\gamma>n$ and $\gamma>1$ here are equivalent to the conditions $\alpha\in(0,1+\frac{2}{n})$ and $\alpha\in(0,n+2)$ in Theorem \ref{thm_affinesphere1}, respectively. We will discuss more about the range of $\gamma$ at the end of this section. Also note that although the domain $\Omega$ here is supposed to be a hyperplane section of the dual cone $C^*$ rather than $C$ itself (see Section \ref{sec_legendre} for details), it is actually easy to see that $\Omega$ satisfies the exterior and the interior sphere conditions if and only if a hyperplane section of $C$ does, so the assumption in  Part \ref{item_analytic2}  coincides with that in the second statement of Theorem \ref{thm_affinesphere1}.


We give the proof of Theorem \ref{thm_affinesphere1analytic} after establishing some lemmata. Although the proof is a simple application of the Monge-Amp\`ere theory, we do not know any reference from which the theorem can be deduced immediately. In fact, most past works on equations of the form \eqref{eqn_affinesphere2} impose strictly convexity and $\C^2$ assumptions on $\pa\Omega$. The improvement to general $\Omega$ in Part \ref{item_analytic1} is achieved by using barrier functions from the following lemma, similarly as in \cite[p.67]{chengyau1}:
\begin{lemma}\label{lemma_simplex}
Let $\Delta\subset\R^n$ be a simplex and $p_0,\cdots,p_n$ be its vertices, so that every $x\in\overline{\Delta}$ can be written uniquely as 
$$
x=t_0(x)p_0+\cdots+t_n(x)p_n,\quad \text{ where } t_0(x),\cdots,t_n(x)\in[0,1],\ t_0(x)+\cdots+t_n(x)=1.
$$ 
Given $\gamma>-n$, consider the function $v\in\C^0(\overline{\Delta})\cap\C^\infty(\Delta)$ defined by
$$
v(x):=-\big(t_0(x)\cdots t_n(x)\big)^{\frac{2}{n+\gamma}},
$$ 
which vanishes on $\pa\Delta$. Then $v$ is convex if and only if $\gamma\geq n$. If $\gamma>n$, then there are constants $C_1,C_2>0$ only depending on $n$, $\gamma$ and the volume of $\Delta$, such that 
$$
C_1 (-v)^{-\gamma}\leq\det\D^2v\leq C_2 (-v)^{-\gamma}
$$
in $\Delta$, and $v$ has infinite slope at every boundary point of $\Delta$.
\end{lemma}
The proof consists of elementary but cumbersome calculations, so we postpone it to the appendix. 
\begin{remark}\label{remark_simplex}
From the expression of $\det\D^2v$ in the proof of Lemma \ref{lemma_simplex}, we see that the equality $\det\D^2v=C(-v)^{-\gamma}$ holds only when $\gamma=n+2$. This gives a solution of Eq.\eqref{eqn_affinesphere2} for $\Omega=\Delta$, $\gamma=n+2$, which corresponds to a hyperbolic affine hypersphere in $\R^{n+1}$ of the form $\{ x_1\cdots x_{n+1}=c\}$ (c.f.\ \cite{calabi}) via Thm.\ \ref{thm_dictionary}.
\end{remark}

The following version of Comparison Principle for Monge-Amp\`ere equations is a special case of \cite[Prop.\@ 2]{chengyau1}. The proof here is adapted from the classical version (see \cite[Thm.\@ 1.4.6]{gutierrez}).
\begin{lemma}\label{lemma_comparison}
Let $\Omega\subset\R^n$ be a bounded convex domain and $u_+,u_-\in\C^0(\pa\Omega)$ be convex functions with boundary values
$\phi_\pm:=u_\pm|_{\pa\Omega}$ such that $\phi_-\leq\phi_+$ and $u_\pm<\env{\phi}_+$ in $\Omega$. Suppose $F(x,t)$ is a positive continuous function on $\{(x,t)\in\Omega\times \R\mid t<\env{\phi}_+(x)\}$ and is non-decreasing in $t$, and $u_+$ (resp.\@ $u_-$) is a generalized supersolution (resp.\@ subsolution) of the Monge-Amp\`ere equation
\begin{equation}\label{eqn_fxu}
\det\D^2u=F(x,u).
\end{equation}
Then we have $u_-\leq u_+$ on $\overline{\Omega}$.
\end{lemma}
\begin{proof}
Suppose by contradiction that $\epsilon:=u_-(x_0)-u_+(x_0)>0$ for some $x_0\in\Omega$. We may take a sufficiently small $\delta>0$ such that $\epsilon-\delta|x-x_0|^2>0$ on $\Omega$. Then
$$E:=\big\{x\in\Omega \,\big|\, u_-(x)-u_+(x)\geq\epsilon-\delta|x-x_0|^2\big\}$$ 
is a compact set with nonempty interior. Since the convex function $u'_-:=u_--\epsilon+\delta|x-x_0|^2$ equals $u_+$ on $\pa E$ and is bounded from below by $u_+$ on $E$, the mass of the Monge-Amp\`ere measure of $u'_-$ on $E$ is no larger than that of $u_+$ (see \eg \cite[Lemma 1.4.1]{gutierrez}). That is, 
\begin{equation}\label{eqn_comparison}
\MA{u_+}(E)\geq\MA{u'_-}(E)=\MA{u_-+\delta|x-x_0|^2}\geq\MA{u_-}(E)+(2\delta)^n\mu(E),
\end{equation}
where $\MA{u}$ denotes the Monge-Amp\`ere measure of $u$, $\mu$ is the Lebesgue measure, and the last inequality follows from the super-additivity of Monge-Amp\`ere measures (see \eg \cite[Lemma 2.9]{figalli}). On the other hand, $u_\pm$ being a generalized super-/sub-solution means $\det\D^2u_+\leq F(x,u_+)$ and $\det\D^2u_-\geq F(x,u_-)$
in the generalized sense, which implies the inequality of measures
$$
F(x,u_+(x))^{-1}\MA{u_+}\leq F(x,u_-(x))^{-1}\MA{u_-}
$$
(where $F(x,u_\pm(x))^{-1}$ are understood as density functions). Since $u_+<u_-$ on $E$ and $F(x,t)$ is non-decreasing in $t$, this implies $\MA{u_+}(E)\leq\MA{u_-}(E)$, 
contradicting \eqref{eqn_comparison}.
\end{proof}

The smoothness assertion in Theorem \ref{thm_affinesphere1analytic} relies on the following regularity result, whose proof is suggested to us by Connor Mooney:
\begin{lemma}\label{lemma_regularity}
Let $\Omega\subset\R^n$ be a bounded convex domain and let $\phi\in\C^0(\pa\Omega)$. Suppose $F(x,t)$ is a smooth function on $\{(x,t)\in\Omega\times\R\mid t<\env{\phi}(x)\}$ and is non-decreasing in $t$, and $u\in\C^0(\overline{\Omega})$ is a convex generalized solution to Eq.\eqref{eqn_fxu} with $u|_{\pa\Omega}=\phi$, such that $u$ has infinite slope at every point of $\pa\Omega$ (this implies $u(x)<\env{\phi}(x)$ for all $x\in\Omega$, see Section \ref{sec_legendre}). Then $u\in\C^\infty(\Omega)$.
\end{lemma}

\begin{proof}
In order to proved that $u$ is smooth around a point $x\in\Omega$, we may assume without loss of generality that $u(x)=0$ and $u\geq0$ on $\overline{\Omega}$. The infinite slope assumption implies that $u>0$ on $\pa\Omega$, hence for sufficiently small $h>0$, the section $S_h(x):=\{y\in\Omega\mid u(y)\leq h\}$ is a convex domain with closure contained in $\Omega$.
	
By the interior regularity result of Cheng-Yau \cite{chengyau1}, if $v\in\C^0(\overline{U})$ is a generalized convex solution to $\det\D^2v=G(x,v)$ in bounded convex domain $U$ with $v|_{\pa U}=0$, where the function $G\in\C^\infty(U\times(-\infty,0])$ is positive, bounded, and non-decreasing in the last variable, then $v\in\C^\infty(U)$\footnote{
The original statement \cite[Thm.\@ 5]{chengyau1} assumes that $U$ satisfies the exterior sphere condition at every boundary point, but this assumption can be removed since it is only used to get paraboloid barriers, while for general $U$ we can use the barriers from Lemma \ref{lemma_simplex} instead. Alternatively, we can use a standard approximation argument (approximate $U$ from the inside by nice subdomains, then apply Pogorelov $\C^2$ estimate and Schauder theory to the solutions on these domains) to deduce the result for general $U$. 
}. 
We obtain the required smoothness of $u$ by applying this to $v=u-h$ on $U=S_h(x)$.
\end{proof}
\begin{remark}\label{remark_regularity}
\begin{enumerate}[label=(\arabic*)]
\item\label{item_regularity1} In the applications below of Lemmas \ref{lemma_comparison} and \ref{lemma_regularity}, the function $F(x,t)$ is either $(-t)^{-\gamma}$ or only a function of $x$. In the former case, a consequence of Lemma \ref{lemma_comparison} that we will often invoke is that if $\Omega_1\subset\Omega_2$ and $w_i\in\C^0(\overline{\Omega}_i)$ is a convex generalized solution to Eq.\eqref{eqn_affinesphere2} for $\Omega=\Omega_i$ ($i=1,2$), then $0\geq w_1\geq w_2$ on $\Omega_1$. In the latter case, Lemma \ref{lemma_comparison} is a special case of the classical Comparison Principle (see \cite[Thm.\@ 1.4.6]{gutierrez}) and the assumption ``$u_\pm<\env{\phi}_+$ in $\Omega$'' is unnecessary.
\item\label{item_regularity2} The assumption that $u$ is continuous up to $\pa\Omega$ is inessential in Lemmas \ref{lemma_comparison} and \ref{lemma_regularity}. Namely, both lemmas are valid for any convex $u\in\C^0(\Omega)$ if we interpret the boundary value $u|_{\pa\Omega}:\pa\Omega\to\R\cup\{+\infty\}$ in the way explained in Section \ref{sec_legendre}. 
Also, the non-decreasing assumption on $F$ is inessential in Lemma \ref{lemma_regularity}, as we can use the interior regularity result of Urbas \cite{urbas} instead of Cheng-Yau's.
\end{enumerate}
\end{remark}

\begin{proof}[Proof of Theorem \ref{thm_affinesphere1analytic} \ref{item_analytic1}]
First suppose that a convex generalized solution $w$ of Eq.\eqref{eqn_affinesphere2} exists. It must be strictly negative in $\Omega$ because otherwise the convexity and the condition $w|_{\pa\Omega}=0$ would imply that $w\equiv0$ in $\Omega$. Therefore, may apply
Lemma \ref{lemma_comparison} to infer that $w$ is unique. Then we can show that $w$ has infinite slope at every $p\in\pa\Omega$ as follows: Let $\Delta\subset\Omega$ be a simplex with a vertex at $p$, $v\in\C^0(\overline{\Delta})\cap\C^\infty(\Delta)$ be the function given by Lemma \ref{lemma_simplex}, and $v_1$ be a constant multiple of $v$ such that $\det\D^2v_1\leq (-v_1)^{-\gamma}$ in $\Delta$. Applying Lemma
\ref{lemma_comparison} to $u_+:=v_1$ and $u_-:=w|_{\overline{\Delta}}$, we get $v_1\geq w$ in $\overline{\Delta}$. Since $v_1$ has infinite slope at $p$ and has the same value as $w$ at $p$, it follows that $w$ also has infinite slope at $p$, as required. By Lemma \ref{lemma_regularity}, it follows that $w\in\C^\infty(\Omega)$.

Thus, we now only need to show the existence of $w$. If $\Omega$ satisfies the exterior sphere condition, this is a special case of \cite[Thm.\@ 1]{chengyau1}. For general $\Omega$, we let $\Omega_1\subset\Omega_2\subset\cdots$ be an exhausting sequence of convex subdomains of $\Omega$ satisfying the exterior sphere condition, and consider the generalized solution $w_i$ on $\Omega_i$. Given $p\in\pa\Omega$, for any simplex $\Delta\subset\R^n$ containing $\Omega$ with $p\in\pa\Delta$, Lemma \ref{lemma_simplex} provides a convex function $w_\Delta\in\C^0(\overline{\Delta})$ vanishing on $\pa\Delta$ and satisfying $\det\D^2w_\Delta\geq(-w_\Delta)^{-\gamma}$, which bounds every $w_i$ from below by Lemma \ref{lemma_comparison} and Remark \ref{remark_regularity}. Therefore, by Arzel\`a-Ascoli Lemma and the fact that the sequence $(w_i)$ is decreasing, the sequence converges uniformly on any compact subset of $\Omega$ to a convex function $w$ on $\Omega$ whose boundary value at $p$ (in the sense of Section \ref{sec_legendre}) is zero. Since $p$ is arbitrary, we actually have $w|_{\pa\Omega}=0$, which also implies that $w\in\C^0(\pa\Omega)$. Now that $w_i$ is a generalized solution and its Monge-Amp\`ere measure weakly converges to that of $w$, we conclude that $w$ is a generalized solution of \eqref{eqn_affinesphere2}, as required.
\end{proof}

Part \ref{item_analytic2} of Theorem \ref{thm_affinesphere1analytic} can be proved by a nearly identical argument, only with the barriers on simplices from Lemma \ref{lemma_simplex} replaced by those on balls from the following lemma, so we omit the details.
\begin{lemma}\label{lemma_ball}
Let $B:=B(0,1)\subset\R^n$ be the unit ball. Then for any $\gamma\geq1$, the radially symmetric function $v\in\C^0(\overline{B})\cap\C^\infty(B)$ given by
$$
v(x):=-(1-|x|^2)^\frac{n+1}{n+\gamma}
$$
is convex, and there are constants $C_1,C_2>0$ only depending on $\gamma$, such that
$$
C_1(-v)^{-\gamma}\leq\det\D^2v\leq C_2(-v)^{-\gamma}
$$
in $B$. Moreover, $v$ has infinite slope at the boundary points if and only if $\gamma>1$.
\end{lemma}
\begin{proof}
For a radially symmetric function $v(x)=f(|x|)$, we have
$$
	\pa_iv(x)=\frac{f'}{|x|}x_i,\quad \pa^2_{ij}v(x)=\frac{f'}{|x|}\left(\delta_{ij}-\frac{x_ix_j}{|x|^2}+\frac{x_ix_j f''}{|x|f'}\right)
$$
($f'$ and $f''$ are evaluated at $|x|$). If $f'(t)>0,f''(t)>0$ for all $0<t<1$, then the matrix
$$
A:=I+\left(\frac{f''}{|x|f'}-\frac{1}{|x|^2}\right)x\, \transp{x}
$$
(where $x$ is viewed as a column vector and $\transp{x}$ is its transpose) is positive definite for $0<|x|<1$ because
$$
\transp{y}Ay=
\begin{cases}
|y|^2&\text{if $y$ is orthogonal to $x$}\\
|x|^2+\left(\frac{f''}{|x|f'}-\frac{1}{|x|^2}\right)|x|^4=\frac{|x|^3f''}{f'}& \text{if $y=x$}.
\end{cases}
$$
The expression of $\pa_{ij}^2v$ can be written as $\D^2v=\frac{f'}{|x|}A$, so $v$ is convex in $B$ this case, and we have
$$
\det\D^2 v(x)=\frac{f'(|x|)^{n-1}f''(|x|)}{|x|^{n-1}}.
$$

Now apply these calculations to $f(t)=-(1-t^2)^\eta$ with $\eta=\frac{n+1}{n+\gamma}$. By a little more calculations, we find that $f'(t),f''(t)>0$ holds for $0<t<1$ if and only if $\eta\leq1$, or equivalently, $\gamma\geq1$, and in this case we have
$$
\det\D^2v(x)=(2\eta)^n\left[1+(1-2\eta)|x|^2\right](-v)^{-\gamma}.
$$
The required inequality follows. The ``Moreover'' statement is elementary to check.
\end{proof}
\begin{remark}
Similarly as Remark \ref{remark_simplex}, the equality $\det\D^2v=C(-v)^{-\gamma}$ holds for the function $v$ in Lemma \ref{lemma_ball} also only when $\gamma=n+2$, and this gives a solution of Eq.\eqref{eqn_affinesphere2} for $\Omega=B$, $\gamma=n+2$, whose corresponding affine hypersphere is the hyperboloid $\{x_1^2+\cdots+x_n^2+1=x_{n+1}^2\}$.
\end{remark}

We also deduce Part \ref{item_analytic3} of Theorem \ref{thm_affinesphere1analytic} from this lemma.
\begin{proof}[Proof of Theorem \ref{thm_affinesphere1analytic} \ref{item_analytic3}]
Fix $0<\gamma\leq1$ and take a ball $B\subset\R^n$ containing $\Omega$ with $\pa B\cap\pa\Omega\neq\emptyset$. By Lemma \ref{lemma_ball}, there is a convex function $v_1\in\C^0(\overline{B})$ which satisfies $\det\D^2v_1\geq (-v_1)^{-\gamma}$ and has finite slopes at boundary points. Any convex generalized solution $u$ to Eq.\eqref{eqn_affinesphere2} is bounded from below by $v_1$ by Lemma \ref{lemma_comparison} while having the same value as $v_1$ at any $p\in\pa B\cap\pa\Omega$, hence has finite slope at $p$. 
\end{proof}

By Parts \ref{item_analytic2} and \ref{item_analytic3}, if $\Omega$ satisfies the exterior and interior sphere conditions, then the condition $\gamma>1$ is necessary and sufficient for the convex solution of Eq.\eqref{eqn_affinesphere2} to have the infinite slope property. However, we do not know whether the condition $\gamma>n$ in Part \ref{item_analytic1} is the optimal sufficient condition for general $\Omega$. An essential problem here is:
\begin{quote}
	\textbf{Question.} When $\Omega\subset\R^n$ is a simplex and $\gamma=n$, does the convex solution $w$ of Eq.\eqref{eqn_affinesphere2} have finite or infinite slope at the vertices of $\Omega$?
\end{quote}
In fact, as $\gamma$ decreases from $>n$ to $0$, we believe that $w$ should start to have finite slope first at the vertices, then at the edges, next at the $2$-facets, and so on (compare Lemma \ref{lemma_triangle} below about the $n=2$ case). So a more general question is: Given $k=0,\cdots,n-1$, what is the exact range of $\gamma$ for $w$ to have infinite slope at every point of the $k$-skeleton of $\Omega$? We leave the study of these questions to future works.

\section{Monge-Amp\`ere problem for affine $(C,k)$-hypersurface}\label{sec_proof2}
 In this section, we prove the following PDE result which implies Theorem \ref{thm_ck1} via Theorem \ref{thm_dictionary}. Note that the equivalence between the assumption on the $C$-regular domain $D=\sepi{\env{\phi}^*}$ in Theorem \ref{thm_ck1} and the assumption on $\phi$ below can be seen by using the rows \circled{1} and \circled{2} in Theorem \ref{thm_dictionary}.
\begin{theorem}[Extended version of Theorem \ref{thm_2}, PDE version of Theorem \ref{thm_ck1}]\label{thm_ck1analytic} 
Let $\gamma>n$, $\Omega\subset\R^n$ be a bounded convex domain and  $w\in\C^0(\overline{\Omega})\cap\C^\infty(\Omega)$ be the convex solution of  Eq.\eqref{eqn_affinesphere2} given by Theorem \ref{thm_affinesphere1analytic}. Suppose $\phi\in\C^0(\pa\Omega)$ satisfies the following condition: for every $p\in\pa\Omega$, there exists an affine function $a:\R^n\to\R$ such that $a(p)=\phi(p)$ and $a\geq\phi$ on $\pa\Omega$. Then for each $\lambda>0$,  a unique convex generalized solution $u$ to 
\begin{equation}\label{eqn_ck1}
\begin{cases}
\det\D^2u=\lambda(-w)^{-\gamma}\ \text{ in }\Omega\\
u|_{\pa\Omega}=\phi
\end{cases}
\end{equation}
exists and has the following properties:
\begin{itemize}
	\item $u$ is smooth in $\Omega$;
	\item $u$ has infinite slope at every point of $\pa\Omega$;
	\item if $u_t$ denotes the solution with parameter $\lambda=e^{-t}$ in Eq.\eqref{eqn_ck1}, then for every fixed $x\in\Omega$, $u_t(x)$ is a strictly increasing concave function in $t$, with value tending to $-\infty$ and $\env{\phi}(x)$ as $t$ tends to $-\infty$ and $+\infty$, respectively.
\end{itemize} 
\end{theorem}
\begin{proof}
The scheme of the proof is the same as Theorem \ref{thm_affinesphere1analytic}. If a generalized solution $u$ of Eq.\eqref{eqn_ck1} exists, then it is unique by Lemma \ref{lemma_comparison} (see also Remark \ref{remark_regularity}), and can be shown to have infinite slope at any $p\in\pa\Omega$ as follows: Let $a$ be an affine function as in the hypothesis and put $w_1:=a+w$. We have $w_1=a\geq \phi=u$ on $\pa\Omega$ and $\det\D^2w_1=\det\D^2u$ in $\Omega$, hence $w_1\geq u$ holds in $\Omega$ by Lemma \ref{lemma_comparison}. Since $w_1$ has infinite slope at $p$ by Theorem \ref{thm_affinesphere1analytic} and has the same value as $u$ at $p$, it follows that $u$ has infinite slope at $p$ as well. Lemma \ref{lemma_regularity} then implies $u\in\C^\infty(\Omega)$.

To show the existence of a generalized solution, we let $\Omega_1\subset\Omega_2\subset\cdots$ be an increasing sequence of strictly convex subdomains of $\Omega$ with $\bigcup_{i}\Omega_i=\Omega$, consider the generalized solution to
$$
\begin{cases}
\det\D^2u=\lambda(-w)^{-\gamma}\ \text{ in }\Omega_i,\\
u|_{\pa\Omega_i}=\env{\phi}|_{\pa\Omega_i},
\end{cases}
$$
which exists because $\int_{\Omega_i}\lambda(-w)^{-\gamma}\dx<+\infty$ (see \eg \cite[Theorem 1.6.2]{gutierrez}), and use the convex function $\env{\phi}+\lambda^{\frac{1}{n}}w\in\C^0(\overline{\Omega})$ as a lower barrier. By the super-additive property of Monge-Amp\`ere measure (see \eg \cite[Lemma 2.9]{figalli}), we have 
$$
\det\D^2(\env{\phi}+\lambda^{\frac{1}{n}}w)\geq \det\D^2\env{\phi}+\det\D^2(\lambda^{\frac{1}{n}}w)=\lambda(-w)^{-\gamma}
$$
in the generalized sense. This allows us to apply Lemma \ref{lemma_comparison} and conclude that $u_i\geq \env{\phi}+\lambda^{\frac{1}{n}}w$ in $\Omega_i$. As in the proof of Theorem \ref{thm_affinesphere1analytic}, using Arzel\`a-Ascoli, we infer that a subsequence of $(u_i)$ converges uniformly on compact subsets to a convex generalized solution $u$ of Eq.\eqref{eqn_ck1} on $\Omega$, as required. An immediate estimate of $u$, which we will use later, is
\begin{equation}\label{eqn_ck3}
\env{\phi}+\lambda^{\frac{1}{n}}w\leq u\leq \env{\phi}.
\end{equation}
In particular, this implies that $u|_{\pa\Omega}=\phi$ in the sense of Section \ref{sec_legendre}. Since $\phi\in\C^0(\pa\Omega)$, we have $u\in\C^0(\overline{\Omega})$.

It remains to show the assertion in the last bullet point. The proof is the same as \cite[\S 8.5]{nie-seppi}. First, using the fact that $\log\det$ is a concave function on the space of positive definite matrices, we get
\begin{align*}
&\log\det \D^2\left(\frac{u_{t_1}+u_{t_2}}{2}\right)=
\log\det\left(\frac{\D^2 u_{t_1}+\D^2u_{t_2}}{2}\right)\\
&\geq \frac{\log\det\D^2u_{t_1}+\log\det\D^2u_{t_2}}{2}=-\frac{t_1+t_2}{2}+\log (-w)^{-\gamma}=\log\det\D^2u_{t_3}
\end{align*}
in $\Omega$ for any $t_1,t_2\in\R$ and $t_3:=\frac{t_1+t_2}{2}$. This implies $\frac{u_{t_1}+u_{t_2}}{2}\leq u_{t_3}$ by Lemma \ref{lemma_comparison}, which means $t\mapsto u_t(x)$ is a concave function for any fixed $x$. Next, 
the required limit $\lim_{t\rightarrow+\infty}u_t(x)=\env{\phi}(x)$ follows from Inequality \eqref{eqn_ck3}, whereas the other limit $\lim_{t\to-\infty}u_t(x)=-\infty$ can be shown by taking a ball $B=B(x,\eps)$ with $\overline{B}\subset\Omega$ and considering the upper barrier
$$
v_t(y)=a\,e^{-\frac{t}{n}}\left(|y-x|^2-\eps^2\right)+b
$$
of $u_t$ in $B$, where the constants $a,b>0$ are chosen to ensure that 
$$
\det\D^2 v_t=(2a)^ne^{-t}\leq e^{-t}(-w)^{-\gamma}=\det\D^2u_t
$$ 
in $B$ and $v_t\geq \env{\phi}\geq u_t$ on $\pa B$ for all $t\in\R$. By Lemma \ref{lemma_comparison} again, the inequality $v_t\geq u_t$ holds in $B$, and in particular at $x$, so the required limit follows. Finally, it follows from Lemma \ref{lemma_comparison} yet again that for any fixed $x\in\Omega$, the function $t\mapsto u_t(x)$ is non-decreasing. If it is not strictly increasing, then the concavity and limiting properties just established would imply that $u_t(x)=\env{\phi}(x)$ for all $t$ larger than some $t_0$. But since  $u_t$ is strictly convex in $\Omega$ (as it is smooth, convex, with $\det\D^2u_t>0$), this is impossible (see Section \ref{sec_legendre}). The proof is completed.
\end{proof}

\begin{remark}\label{remark_boundary}
\begin{enumerate}[label=(\arabic*)]
\item Since the assumption $\gamma>n$ is inherited from Part \ref{item_analytic1} of Theorem \ref{thm_affinesphere1analytic},  by using Part  \ref{item_analytic2} instead, we may relax it to $\gamma>1$ if $\Omega$ satisfies both the exterior and interior sphere conditions.
\item
If $\Omega$ is a bounded convex domain with $\C^2$ boundary and satisfying the exterior sphere condition, then a sufficient condition for $\phi$ to fulfill the assumption in Theorem \ref{thm_ck1analytic} is $\phi\in\C^2(\pa\Omega)$. Indeed, fix a point $p$ and write $\pa\Omega$ locally as a graph of the form $\{x+f(x)v\mid x\in U\}$, where $U\subset H$ is an open subset, $H$ is the tangent hyperplane of $\pa\Omega$ through $p$, $v$ is its inward unit normal vector, and $f\in C^2(U)$ is a convex function. Then the exterior sphere condition at $p$ is equivalent to the condition that the Hessian of $f$ is positive definite at $p$. Up to adding an affine function, we can suppose that $\varphi$ has a critical point at $p$; now one easily sees that the affine function mapping $x+tv$ to $x+c t v$ is larger than $\varphi$ on $\pa\Omega$ if $c>0$ is sufficiently large.
\item
As one can see from the proof, the existence and uniqueness of the convex generalized solution $u$ to Eq.\eqref{eqn_ck1} hold independently of the assumption on $\phi$. The assumption is only used to ensure the regularity and the infinite slope property.
\item 
When $\Omega$ is the unit ball in $\R^n$ with $n\geq3$, the examples constructed by Bonsante-Fillastre \cite[\S 3.7]{bonsante-fillastre} correspond to certain $\phi\in\C^0(\pa\Omega)$ such that the convex generalized solution $u$ of Eq.\eqref{eqn_ck1} is not smooth or strictly convex. In fact, $u$ restricts to affine functions on many line segments in $\Omega$ joining boundary points, and coincides with $\env{\phi}$ on these segments.   
\end{enumerate}
\end{remark}

\section{Monge-Amp\`ere problem for affine $(C,k)$-surface ($n=2$ case).}\label{sec_proof3}
In this final section, we generalize Theorem \ref{thm_ck1analytic} along the line of our previous work \cite{nie-seppi} and obtain Theorem \ref{thm_ck2analytic} below, which implies Theorem \ref{thm_ck2} ($=$ Theorem \ref{thm_2} in the introduction) via Theorem \ref{thm_dictionary}.

Recall that the main novelty of \cite{nie-seppi} on the PDE aspect is the study of Eq.\eqref{eqn_ck1} when $n=2$ and $\phi$ is merely an $\R\cup\{+\infty\}$-valued lower semicontinuous function. The simplest nontrivial example is 
\begin{equation}\label{eqn_trianglefunction}
\phi(x)=
\begin{cases}
0&\text{ if }x\in\{p_1,p_2,p_3\},\\
+\infty& \text{ otherwise},
\end{cases}
\end{equation}
where $p_1,p_2,p_3\in\pa\Omega$. Another typical example is $\phi=0$ throughout $\pa\Omega$ except at a single point $p$, with $\phi(p)=-1$. When $\Omega$ is the unit disk, the solution of Eq.\eqref{eqn_ck1} for the latter example is computed in \cite{bs}.

When $\phi=+\infty$ on a part of $\pa\Omega$, we interpret Eq.\eqref{eqn_ck1} as the problem of seeking a lower semicontinuous convex function $u$ on $\R^2$ with $u|_{\pa\Omega}=\phi$, such that $\dom{u}$ is contained in $\overline{\Omega}$ and has nonempty interior $U$, with the Monge-Amp\`ere equation satisfied in $U$. We shall add to this problem the extra constraint that $u$ has infinite slope at every boundary point of $U$ within $\Omega$, otherwise the solution would not be unique for trivial reason: for example, if $\phi$ is as in \eqref{eqn_trianglefunction}, then for any closed convex set $E\subset\overline{\Omega}$ with $E\cap\pa\Omega=\{p_1,p_2,p_3\}$, modifying the values of the function $w$ outside of $E$ into $+\infty$  would give a solution.

The main result below asserts that Eq.\eqref{eqn_ck1}, interpreted in this way, has a unique solution. Here the spaces of functions $\LC(\R^2)$ and $\LC(\pa\Omega)$ were introduced in Section \ref{sec_legendre}. Also recall that the essential domain $\dom{\env{\phi}}$ of the convex envelope $\env{\phi}$ of any $\phi\in\LC(\pa\Omega)$ coincides with the convex hull of $\dom{\phi}$.
\begin{theorem}[Extended version of Theorem \ref{thm_1}, PDE version of Theorem \ref{thm_ck2}]\label{thm_ck2analytic}
Let $\gamma>2$, $\Omega\subset\R^2$ be a bounded convex domain, and  $w\in\C^0(\overline{\Omega})\cap\C^\infty(\Omega)$ be the convex solution of  Eq.\eqref{eqn_affinesphere2} given by Theorem \ref{thm_affinesphere1analytic}. Let $\phi\in\LC(\pa\Omega)$ be such that $\dom{\env{\phi}}$ has nonempty interior. Then for any $\lambda>0$, there exists a unique $u\in\LC(\R^2)$ satisfying
	\begin{equation}\label{eqn_ck2}
	\begin{cases}
	U:=\interior\dom{u}\neq\emptyset,\ U\subset\Omega,\\
	\det\D^2u=\lambda(-w)^{-\gamma}\ \text{ in $U$},\\
	u|_{\pa\Omega}=\phi,\\
	\text{$u$ has infinite slope at every point of $\pa U\cap\Omega$,}
	\end{cases}
	\end{equation}
in the generalized sense, and it has the following properties:
	\begin{enumerate}[label=(\alph*)]
		\item\label{item_ck41} $u$ is smooth in $U$;
		\item\label{item_ck42} there exists a convex function $f\in\C^0(\overline{U})$ with $f|_{\pa U}=0$ such that $\env{\phi}+f\leq u\leq\env{\phi}$ on $\overline{U}$;
		\item\label{item_ck44} if $u_t$ denotes the solution with parameter $\lambda=e^{-t}$ in Eq.\eqref{eqn_ck2}, then for every fixed $x\in U$, $u_t(x)$ is a strictly increasing concave function in $t$, with value tending to $-\infty$ and $\env{\phi}(x)$ as $t$ tends to $-\infty$ and $+\infty$, respectively.
	\end{enumerate}	
Moreover, if $\gamma\geq4$ and $\Omega$ satisfies the exterior circle condition at $p\in\pa U\cap\pa\Omega$, then $u$ has infinite slope at $p$.
\end{theorem}
In particular, if $\Omega$ satisfies the exterior circle condition at every boundary point, then the resulting solutions belong to the space $\SC_0(\Omega)$ defined in Section \ref{sec_legendre}. Therefore, using Theorem \ref{thm_dictionary}, one readily checks that this theorem implies the geometric result, Theorem \ref{thm_ck2}.

Before giving the proof of Theorem \ref{thm_ck2analytic}, we first note that for a general $u\in\LC(\R^2)$ with $u|_{\pa\Omega}=\phi$, since $u\leq\env{\phi}$ on $\R^2$, the convex set $\dom{u}$ contains $\dom{\env{\phi}}$, but the two sets might not coincide. Even if they coincide, the values of $u$ and $\env{\phi}$ on the boundary of the set might not be the same. However, we do have $\dom{u}=\dom{\env{\phi}}$ and $u|_{\pa U}=\env{\phi}|_{\pa U}$ for $u$ satisfying \eqref{eqn_ck2}. The proof of this relies on and the following Generalized Comparison Principle from \cite{nie-seppi}:

\begin{lemma}[{\cite[Lemma 6.4]{nie-seppi}}]\label{lemma_generalizedcomparison}
	Let $U\subset\R^n$ be a bounded convex domain,
	$u_+:\overline{U}\rightarrow\Rp$ be a lower semicontinuous convex function taking finite values in $U$ and $u_-:\overline{U}\rightarrow\R$ be a continuous convex function such that
	\begin{itemize}
		\item $\det\D^2u_+\leq \det\D^2u_-$ in $U$ in the generalized sense;
		\item $u_+(x)\geq u_-(x)$ for every $x\in\pa U$ where $u_+$ has finite slope, as well as every $x\in\pa U$ where $u_-$ has infinite slope.
	\end{itemize}
	Then we have $u_+\geq u_-$ throughout $\overline{U}$.
\end{lemma}

The infinite slope property in Theorem \ref{thm_ck2analytic} will be treated by using the following results:
\begin{lemma}[{\cite[Lemmas 8.3 and 8.4]{nie-seppi}}]\label{lemma_triangle}
Let $\Delta\subset\R^2$ be a triangle, $p\in\pa\Delta$ be a vertex, $I\subset\pa\Delta$ be an edge, and $u\in\C^0(\overline{\Delta})$ be a convex function. Then the following statements hold.
\begin{enumerate}[label=(\arabic*)]
	\item\label{item_triangle1}
	If there are constants $c>0$ and $\beta>-2$ such that
	$$
	\begin{cases}
	\det\D^2 u(x)\leq c\,|x-p|^\beta \text{ for }x\in \Delta\\
	u|_{\pa\Delta}=0
	\end{cases}
	$$
    in the generalized sense, then $u$ has finite slope at $p$.
	\item\label{item_triangle2}
	If there is a constant $c>0$ such that
	$$
	\det\D^2 u(x)\geq c\,|x-p|^{-2}\ \text{ for }x\in \Delta
	$$
	 in the generalized sense, then $u$ has infinite slope at $p$.
	\item\label{item_triangle3}
	If there is a constant $c>0$ such that
	$$
    \begin{cases}
\det\D^2 u\geq c\ \text{ in }\Delta\\
u|_I=0
\end{cases}
     $$
	in the generalized sense, then $u$ has infinite slope at every interior point of $I$.
\end{enumerate}
\end{lemma}

We now present the proof of Theorem \ref{thm_ck2analytic}. Since it is very similar to the $\gamma=4$ case treated in \cite[Theorem A']{nie-seppi} except a few modifications, we do not give all the details.
\begin{proof}[Proof of Theorem \ref{thm_ck2analytic}]
First, we fix a solution $u\in\LC(\R^2)$ of Eq.\eqref{eqn_ck2} with parameter $\lambda=1$ and prove properties \ref{item_ck41} and \ref{item_ck42}. The proof easily adapts to general $\lambda>0$.

Since $(-w)^{\gamma}$ is a positive smooth function in $\Omega$, the regularity theory of Monge-Amp\`ere equations in two variables (see \eg \cite[Thm.\@ 6.7]{nie-seppi}) immediately implies property \ref{item_ck41}.
	
We proceed to show property \ref{item_ck42} when $f=w_0\in\C^0(\overline{U})\cap\C^\infty(U)$ is the convex solution to 
\begin{equation}\label{eqn_proofck41}
\begin{cases}
\det\D^2w_0=(-w_0)^{-\gamma}\ \text{ in }U,\\
w_0|_{\pa U}=0,
\end{cases}
\end{equation}
(whose unique-existence is given by Theorem \ref{thm_affinesphere1analytic}). Namely, we shall show 
\begin{equation}\label{eqn_proofck42}
u\geq\env{\phi}+w_0\ \text{ on }\overline{U}.
\end{equation}

For the sake of clearness, let us first show the weaker inequality $u\geq\env{\phi}+w$ on $\overline{U}$. Since $\env{\phi}$ is the pointwise supremum of all affine functions $a:\R^2\to\R$ satisfying $a|_{\pa\Omega}\leq\phi$, it suffices to show $u\geq a+w$ for any such $a$. But this follows from Lemma \ref{lemma_generalizedcomparison} because on one hand, the Monge-Amp\`ere measures of $u$ and $a+w$ coincide; on the other hand, a point $x\in\pa U$ where either $u$ has finite slope or $a+w$ has infinite slope can only lie on $\pa\Omega$, where we have $u(x)\geq a(x)=a(x)+w(x)$ by assumptions. 

This argument does not work immediately if we replace $w$ by $w_0$: although we still have the correct comparison of Monge-Amp\`ere measures required by Lemma \ref{lemma_generalizedcomparison}, namely 
$$
\det\D^2u=(-w)^{-\gamma}\geq (-w_0)^{-\gamma}=\det\D^2(a+w_0)
$$ 
(this follows from Lemma \ref{lemma_comparison}), we lack the comparison of boundary values, as $a+w_0$ now has infinite slope on the whole $\pa U$ rather than just on $\pa U\cap\pa\Omega$. Nevertheless, it can be adjusted as follows. Let $U_\delta:=\{x\in\Omega\mid d(x,U)<\delta\}$ be the intersection of the $\delta$-neighborhood of $U$ with $\Omega$, and $w_\delta\in\C^0(\overline{U}_\delta)\cap\C^\infty(U_\delta)$ be the solution of Eq.\eqref{eqn_proofck41} with $U$ replaced by $U_\delta$. Then the argument is valid for $w_\delta$ and yields
\begin{equation}\label{eqn_proofck43}
u\geq a+w_\delta\ \text{ on }\overline{U}.
\end{equation}
As $\delta$ tends to $0$, the restriction of $w_\delta$ to $U$ increases and converges uniformly to $w_0$ (the proof of this is given in \cite[Prop.\@ 7.8]{nie-seppi} for the $\gamma=n+2$ case, which can be generalized to any $\gamma> n$ by using the barrier functions from Lemma \ref{lemma_simplex}). Therefore, by taking the pointwise supremum of \eqref{eqn_proofck43} for all affine function $a$ with $a|_{\pa\Omega}\leq\phi$ and all $\delta>0$, we obtain the required inequality \eqref{eqn_proofck42} and finish the proof of \ref{item_ck42}.

Next, we show the unique-existence of solution $u\in\LC(\R^2)$ to Eq.\eqref{eqn_ck2}. By property \ref{item_ck42}, if we let $U$ denote the interior of $\dom{\env{\phi}}$, then a solution $u\in\LC(\R^2)$ of \eqref{eqn_ck2} is equivalent to a convex generalized solution $\widetilde{u}\in\C^0(U)$ of
\begin{equation}\label{eqn_tildeu}
\begin{cases}
\det\D^2\widetilde{u}=(-w)^{-\gamma}\ \text{ in }U,\\
\widetilde{u}|_{\pa U}=\env{\phi}|_{\pa U},\\
\text{$\widetilde{u}$ has infinite slope at every point of $\pa U\cap\Omega$},
\end{cases}
\end{equation}
(we obtain $u$ by setting the value of $\widetilde{u}$ outside $\overline{U}$ to be $+\infty$). If the last slope condition in Eq.\eqref{eqn_tildeu} is removed, one easily shows the unique-existence of $\widetilde{u}$ by a standard approximation argument, as in the proofs of Theorems \ref{thm_affinesphere1analytic} and \ref{thm_ck1analytic}, using the function $\env{\phi}+w_0$ as a lower barrier. But such a $\widetilde{u}$ automatically has infinite slope at any $q\in\pa U\cap\Omega$ because of Lemma \ref{lemma_triangle} \ref{item_triangle3} and the fact that $q$ lies on a line segment in $\pa U$ (since $\dom{\env{\phi}}$ is the convex envelope of points on $\pa\Omega$). This shows the unique-solvability of \eqref{eqn_tildeu}, and hence that of \eqref{eqn_ck2}.

Now it only remains to show \ref{item_ck44} and the last ``Moreover'' statement. Property \ref{item_ck44} has already been proven in Theorem \ref{thm_ck1analytic} in a slightly different setting. The proof still works in the current setting, so we omit the details. The idea of proof for the last statement is to estimate $w$ near $p$ using the function from Lemma \ref{lemma_ball} as lower barrier, then apply Lemma \ref{lemma_triangle} \ref{item_triangle2}. More precisely, let $B\subset\R^2$ be a disk containing $\Omega$ such that $\pa B$ passes through $p$. Assuming without loss of generality that
$B=B(0,1)$, we use Lemma \ref{lemma_ball} to get a convex function $v\in\C^0(\overline{B})$ of the form $v=-C(1-|x|^2)^\frac{3}{2+\gamma}$ satisfying $\det\D^2v\geq (-v)^{-\gamma}$, which bounds $w$ from below by Lemma \ref{lemma_comparison}. Given a triangle $\Delta\subset U$ with a vertex at $p$, we then have the following estimate for the right-hand side of the Monge-Amp\`ere equation for $u$ on $\Delta$: 
$$
(-w)^{-\gamma}\geq (-v)^{-\gamma}=C^{-\gamma}(1-|x|^2)^{-\frac{3\gamma}{2+\gamma}}\geq C_1|x-p|^{-\frac{3\gamma}{2+\gamma}}\geq C_2|x-p|^{-2}\quad (\text{for all }x\in\Delta),
$$ 
where $C_1,C_2>0$ only depends on $\gamma$ and $\Delta$, the last inequality uses the assumption $\gamma\geq4$, and the second-to-last inequality is because we have $c|x-p|\leq 1-|x|^2\leq c^{-1}|x-p|$ for all $x\in\Delta$ and some $0<c<1$ only depending on $\Delta$. Therefore, we obtain the required statement by applying Lemma \ref{lemma_triangle} \ref{item_triangle2} to $w|_{\overline{\Delta}}$.
\end{proof}

Finally, the following result corresponds to Proposition \ref{prop_alpha1} and explains the necessity of the condition $\gamma\geq4$ in the last statement of Theorem \ref{thm_ck2analytic}:
\begin{proposition}[PDE version of Prop.\ \ref{prop_alpha1}]\label{prop_circle}
Under the hypotheses of Theorem \ref{thm_ck2analytic}, further assume that $\gamma<4$, $\Omega=B:=B(0,1)$ is the unit disk, and $\phi$ is the function given by \eqref{eqn_trianglefunction}, namely $\phi=0$ at $p_1,p_2,p_3\in\pa B$ and $\phi=+\infty$ everywhere else. 
Then the solution $u\in\LC(\R^2)$ to Eq.\eqref{eqn_ck2} has finite slope at $p_1$, $p_2$ and $p_3$.
\end{proposition}
\begin{proof}
Let us show that $u$ has finite slope at $p_1$. Let $\Delta$ denote the triangle with vertices $p_1$, $p_2$ and $p_3$, so that $\dom{u}=\dom{\env{\phi}}=\overline{\Delta}$ by property \ref{item_ck42} in Theorem \ref{thm_ck2analytic}.
Lemma \ref{lemma_ball} gives a convex function $v\in\C^0(\overline{B})$ of the form $v(x)=-C(1-|x|^2)^{\frac{3}{2+\gamma}}$ satisfying $\det\D^2v\leq (-v)^{-\gamma}$, which bounds $w$ from above by Lemma \ref{lemma_comparison}. Similarly as in the last paragraph of the previous proof, we get the following estimate for the right-hand side of the Monge-Amp\`ere equation for $u$:
$$
(-w)^{-\gamma}\leq (-v)^{-\gamma}=C^{-\gamma}(1-|x|^2)^{-\frac{3\gamma}{2+\gamma}}\leq C_1|x-p_1|^{-\frac{3\gamma}{2+\gamma}}\quad (\text{for all }x\in\Delta).
$$ 
Since we have $-\frac{3\gamma}{2+\gamma}>-2$ by the assumption $\gamma<4$ and $u|_{\pa\Delta}=0$ by Theorem \ref{thm_ck2analytic}, we can apply Lemma \ref{lemma_triangle} \ref{item_triangle1} and conclude that $u$ has finite slope at $p_1$.
\end{proof}

\appendix
\section{Proof of Lemma \ref{lemma_simplex}}
	By applying an affine transformation to the simplex $\Delta$, which changes the Monge-Amp\`ere measure of any convex function by a constant factor only depending on the Jacobian of the affine transformation, we may assume that the vertices of $\Delta$ are 
	$$
	p_0=0,\ p_1=(1,0,\cdots,0),\ \cdots,\ p_n=(0,\cdots,0,1),
	$$
	so that we have $t_0(x)=1-x_1-\cdots-x_n$ and $t_i(x)=x_i$. This reduces Lemma \ref{lemma_simplex} to:
\begin{lemma}
Let $\Delta:=\big\{x\in\R^n\,\big|\,x_1,\cdots,x_n>0\, x_1+\cdots+x_n<1\big\}$. Given $\gamma>-n$, let $v\in\C^0(\overline{\Delta})$ be defined by
	$$
	v(x):=-\big(x_0 x_1\cdots x_n\big)^{\frac{2}{n+\gamma}},
	$$ 
	where given any $x\in\overline{\Delta}$, we denote $x_0:=1-x_1-\cdots-x_n$. Then $v$ is convex if and only if $\gamma\geq n$. If $\gamma>n$, then there are constants $C_1,C_2>0$ only depending on $n$ and $\gamma$, such that 
	\begin{equation}\label{eqn_app}
	C_1 (-v)^{-\gamma}\leq\det\D^2v\leq C_2 (-v)^{-\gamma}
	\end{equation}
	in $\Delta$, and $v$ has infinite slope at every boundary point of $\Delta$.
\end{lemma}
\begin{proof}
The last infinite slope property when $\gamma>n$ is elementary to check using the expression of $v$. To prove the other assertions, let us first compute $\det\D^2v$. Fix $\gamma>-n$ and put
	$$
	\eta:=\tfrac{2}{n+\gamma}~,\quad V(x):=\eta\big(\log x_0+\log x_1+\cdots+\log x_n\big).
	$$
	Then we have $v=-e^V$ and 
	$$
	\pa^2_{ij}v=v\big(\pa^2_{ij}V+\pa_iV\pa_jV\big)=:(-v)\big(A_{ij}-B_iB_j\big),
	$$
	where we set
	$$
	B_i:=\pa_i V=\eta\left(\frac{1}{x_i}-\frac{1}{x_0}\right),\quad A_{ij}:=- \pa^2_{ij}V=\eta\left(\frac{\delta_{ij}}{x_i^2}+\frac{1}{x_0^2}\right).
	$$
	
We have the following equalities for the determinant of a symmetric block matrix:
		$$
		\det\begin{pmatrix}
		A&B\\
		\transp{B}&C
		\end{pmatrix}
		=\det(A)\det(C-\transp{B}A^{-1}B)=\det(C)\det(A-BC^{-1}\transp{B})
		$$
(the first equality can be shown by eliminating the upper-right block, and the second by eliminating the lower-left one). As a consequence, when $C$ is the $1\times1$ identity, we get
		\begin{equation}\label{eqn_app1}
		\det(A-B\,\transp{B})=(1-\transp{B}A^{-1}B)\det(A).
		\end{equation}
In order to apply this formula to the above-defined $A$ and $B$, we need to compute $\transp{B}A^{-1}B$ and $\det(A)$. To this end, note that $A$ can be written as
		$$
		A=\frac{ \eta}{x_0^2}\left[
		\begin{pmatrix}
		(x_0/x_1)^2&&\\
		&\ddots&\\
		&&(x_0/x_n)^2
		\end{pmatrix}
		+
		\begin{pmatrix}
		1\\
		\vdots\\
		1
		\end{pmatrix}
		(1,\cdots,1)
		\right]=:\frac{ \eta}{x_0^2}(\Lambda+E),
		$$
It follows that $A^{-1}$ can be calculated through series expansion as
		\begin{align*}
	A^{-1}&=\frac{x_0^2}{ \eta}\big(\Lambda^{-1}-\Lambda^{-1}E\Lambda^{-1}+\Lambda^{-1}E\Lambda^{-1}E\Lambda^{-1}-\cdots\big)=\frac{x_0^2}{ \eta}\sum_{k\geq0}(-1)^k(\Lambda^{-1}E)^k\Lambda^{-1}.
	\end{align*}
When $k\geq1$, we have
$$
(\Lambda^{-1}E)^k\Lambda^{-1}=\Lambda^{-1}	
\begin{pmatrix}
1\\
\vdots\\
1
\end{pmatrix}
\left[
(1,\cdots,1)
\Lambda^{-1}
	\begin{pmatrix}
1\\
\vdots\\
1
\end{pmatrix}
\right]^{k-1}
\hspace{-0.2cm}
(1,\cdots,1)\Lambda^{-1}
=\left(\frac{|x|^2}{x_0^2}\right)^{k-1}\Lambda^{-1}E\Lambda^{-1},
$$	
where $|x|^2:=x_1^2+\cdots+x_n^2$. Therefore, we obtain
$$
	A^{-1}=\frac{x_0^2}{ \eta}\left(\Lambda^{-1}-\frac{x_0^2}{x_0^2+|x|^2}\Lambda^{-1}E\Lambda^{-1}\right)=\frac{1}{ \eta}\begin{pmatrix}
		x_1^2&&\\
		&\ddots&\\
		&&x_n^2
		\end{pmatrix}
		-\frac{1}{ \eta(x_0^2+|x|^2)}
		\begin{pmatrix}
		x_1^2\\
		\vdots\\
		x_n^2
		\end{pmatrix}
		(x_1^2,\cdots,x_n^2),
$$
$$
		\transp{B}A^{-1}B= \eta\left\{\sum_{i=1}^n\left(1-\frac{x_i}{x_0}\right)^2-\frac{1}{x_0^2+|x|^2}\left[\sum_{i=1}^n\left(1-\frac{x_i}{x_0}\right)x_i\right]^2\right\}~.
$$
The two sums in the expression of $\transp{B}A^{-1}B$ can be simplified as follows:
		$$
		\sum_{i=1}^n\left(1-\frac{x_i}{x_0}\right)^2=\sum_{i=1}^n\left(1+\frac{x_i^2}{x_0^2}-\frac{2x_i}{x_0}\right)=n+\frac{|x|^2}{x_0^2}-\frac{2(1-x_0)}{x_0}=n+2+\frac{|x|^2}{x_0^2}-\frac{2}{x_0}~,
		$$
		$$
		\sum_{i=1}^n\left(1-\frac{x_i}{x_0}\right)x_i=(1-x_0)-\frac{|x|^2}{x_0}=1-\frac{x_0^2+|x|^2}{x_0}~.
		$$
It follows that
$$
		\transp{B}A^{-1}B= \eta\left\{n+2+\frac{|x|^2}{x_0^2}-\frac{2}{x_0}-\frac{1}{x_0^2+|x|^2}\left(1-\frac{x_0^2+|x|^2}{x_0}\right)^2\right\}=\eta\left(n+1-\frac{1}{x_0^2+|x|^2}\right)
$$

On the other hand, noting  that \eqref{eqn_app1} still holds if both ``$-$'' signs are replaced by ``$+$'', we get
$$
\det(A)=\left(\frac{\eta}{x_0^2}\right)^n\left[1+\left(\frac{x_1}{x_0}\right)^2+\cdots+\left(\frac{x_1}{x_0}\right)^2\right]\left(\frac{x_0^n}{x_1\cdots x_n}\right)^2=\frac{ \eta^n(x_0^2+|x|^2)}{(x_0x_1\cdots x_n)^2}~.
$$
Putting everything together, we obtained the expression of $\det\D^2v$ as
	\begin{align*}
	\det\D^2 v&=(-v)^n\det(A-B\,\transp{B})=(-v)^n(1-\transp{B}A^{-1}B)\det(A)\\
	&=(x_0x_1\cdots x_n)^{ \eta n}\left(1- \eta(n+1)+\frac{ \eta}{x_0^2+|x|^2}\right)\frac{ \eta^n(x_0^2+|x|^2)}{(x_0x_1\cdots x_n)^2}\\
	&=\eta^{n+1}\left[1-\big(n+1-\eta^{-1}\big)\big(x_0^2+|x|^2\big)\right](-v)^{-\gamma}.
	\end{align*}
	
The condition $\gamma>n$ is equivalent to $n+1-\eta^{-1}<1$. On the other hand, it is elementary to check that $x_0^2+|x|^2\leq1$ for all $x\in\overline{\Delta}$, with equality exactly when $x$ is a vertex. Thus, we obtain the required inequality \eqref{eqn_app} when $\gamma>n$. We may also conclude that if $\gamma\geq n$, then $\det\D^2v>0$ throughout $\Delta$, otherwise we have $\det\D^2v<0$ near a vertex. In particular, $v$ is not convex if $-n<\gamma<n$.

Now it only remains to be shown that $v$ is indeed convex when $\gamma\geq n$. We shall fix such a $\gamma$ and show that $\D^2v>0$ in $\Delta$. Since we already know $\det\D^2v>0$ in $\Delta$, if $\D^2v(x)$ is not positive definite for $x\in\Delta$, then it has at least two negative eigenvalues. But by the above computations, $\D^2v(x)$ can be written as
$\widetilde{A}-\widetilde{B}\,\transp{\widetilde{B}}$ for the positive definite symmetric matrix $\widetilde{A}=-v(x)A(x)$ and column vector $\widetilde{B}=-v(x)B(x)$ (note that $v(x)<0$). Hence it admits at most one negative eigenvalue, because $\transp{y}(\widetilde{A}-\widetilde{B}\,\transp{\widetilde{B}})y=\transp{y}\widetilde{A}y>0$ for any $y$ orthogonal to $\widetilde{B}$. This contradiction shows $\D^2v>0$ and completes the proof.
\end{proof}

\bibliographystyle{siam} 
\bibliography{deformation}
\end{document}